\documentclass[11pt]{article}
\usepackage{amstex}

\newcount\sectiondepthcount
\newcommand{\sectiondepth}[1]
  {\sectiondepthcount=#1
   \ifnum\sectiondepthcount=3
      \numberwithin{equation}{subsubsection}
   \fi
   \ifnum\sectiondepthcount=2
      \numberwithin{equation}{subsection}
   \fi}
\sectiondepth{2}

\renewcommand{\theequation}{%
  \ifnum\sectiondepthcount=3%

\arabic{section}.\arabic{subsection}.\arabic{subsubsection}.\arabic{equation}%
    \fi
 \ifnum\sectiondepthcount=2
    \arabic{section}.\arabic{subsection}.\arabic{equation}\fi}

\renewcommand{\thesection}{\arabic{section}}

\renewcommand{\section}[1]
  {\vskip-\lastskip\penalty-1000\vskip\bigskipamount\refstepcounter{section}
   \noindent{\bf\thesection.\ #1} \vskip\bigskipamount\par}

\renewcommand{\thesubsection}{\arabic{section}.\arabic{subsection}}

\renewcommand{\subsection}[1]
 {\vskip-\lastskip\penalty-500\vskip\bigskipamount\refstepcounter{subsection}
  \ifnum\sectiondepthcount=3{\bf\thesubsection.\ }{\bf
#1}\vskip\bigskipamount\par\fi
  \ifnum\sectiondepthcount=2{\bf(\thesubsection)\ }{\em #1}\fi}

\renewcommand{\thesubsubsection}%
    {\arabic{section}.\arabic{subsection}.\arabic{subsubsection}}

\renewcommand{\subsubsection}[1]
  {\vskip-\lastskip\vskip\bigskipamount\refstepcounter{subsubsection}
   {\bf (\thesubsubsection)\ }}

\newenvironment{claim}[1]
  {\vskip-\lastskip\penalty-200\vskip\bigskipamount%
   \ifnum\sectiondepthcount=3
      \refstepcounter{subsubsection}%
      {\bf (\arabic{section}.\arabic{subsection}.\arabic{subsubsection}) }%
      \sc{#1}--- \it%
   \fi
   \ifnum\sectiondepthcount=2
      \refstepcounter{subsection}%
      {\bf (\arabic{section}.\arabic{subsection}) }%
      \sc{#1}--- \it%
   \fi}
  {\vskip-\lastskip\penalty+200\vskip\bigskipamount\rm}

\newenvironment{claim*}[1]{\vskip-\lastskip\penalty-200\vskip\bigskipamount%
   \sc{#1}--- \it}{\vskip-\lastskip\penalty+200\vskip\bigskipamount\rm}

\newenvironment{claimwem}[1]
  {\vskip-\lastskip\penalty-200\vskip\bigskipamount%
   \ifnum\sectiondepthcount=3
      \refstepcounter{subsubsection}%
      {\bf (\arabic{section}.\arabic{subsection}.\arabic{subsubsection}) }%
      \sc{#1}--- \rm%
   \fi
   \ifnum\sectiondepthcount=2
      \refstepcounter{subsection}%
      {\bf (\arabic{section}.\arabic{subsection}) }%
      \sc{#1}--- \rm%
   \fi}
  {\vskip-\lastskip\penalty+200\vskip\bigskipamount\rm}

\newenvironment{claimwem*}[1]{\vskip-\lastskip\penalty-500\vskip\bigskipamount%
   \sc{#1}---\rm}{\rm\penalty+500\par\vskip\bigskipamount\ignorespaces}

\newcommand{\theoremname}{Theorem.}
\newcommand{\corollaryname}{Corollary.}
\newcommand{\lemmaname}{Lemma.}
\newcommand{\propositionname}{Proposition.}
\newcommand{\conjecturename}{Conjecture.}
\newcommand{\definitionname}{Definition.}
\newcommand{\examplename}{Example.}
\newcommand{\remarkname}{Remark.}
\newcommand{\pfname}{Proof.}

\renewenvironment{th}{\begin{claim}{\theoremname}}{\end{claim}}
\newenvironment{cor}{\begin{claim}{\corollaryname}}{\end{claim}}
\newenvironment{lem}{\begin{claim}{\lemmaname}}{\end{claim}}
\newenvironment{prop}{\begin{claim}{\propositionname}}{\end{claim}}
\newenvironment{defi}{\begin{claim}{\definitionname}}{\end{claim}}

\newenvironment{rem}{\begin{claimwem}{\remarkname}}{\end{claimwem}}

\newenvironment{th*}[1]{\begin{claim*}{\theoremname #1}}{\end{claim*}}
\newenvironment{cor*}[1]{\begin{claim*}{\corollaryname #1}}{\end{claim*}}
\newenvironment{lem*}[1]{\begin{claim*}{\lemmaname #1}}{\end{claim*}}
\newenvironment{prop*}[1]{\begin{claim*}{\propositionname #1}}{\end{claim*}}
\newenvironment{defi*}[1]{\begin{claim*}{\definitionname #1}}{\end{claim*}}
\newenvironment{conj*}[1]{\begin{claim*}{\conjecturename #1}}{\end{claim*}}
\newenvironment{rem*}{\begin{claimwem*}{\remarkname}}{\end{claimwem*}}
\newenvironment{ex*}{\begin{claimwem*}{\examplename}}{\end{claimwem*}}

\newenvironment{proof}{\vskip-\lastskip\vskip\medskipamount{\it\pfname}}%
                       {$\square$\vskip\medskipamount\par}

\def\thebibliography#1{\vskip-\lastskip\penalty-1000\vskip\bigskipamount%
\centerline{{\bf\refname}}\bigskip
  \list{[\arabic{enumi}]}{\settowidth\labelwidth{[#1]}%
  \leftmargin\labelwidth \advance\leftmargin\labelsep
  \usecounter{enumi}}}

\newskip\tmpskip

\mathcode`A="7041 \mathcode`B="7042 \mathcode`C="7043 \mathcode`D="7044
\mathcode`E="7045 \mathcode`F="7046 \mathcode`G="7047 \mathcode`H="7048
\mathcode`I="7049 \mathcode`J="704A \mathcode`K="704B \mathcode`L="704C
\mathcode`M="704D \mathcode`N="704E \mathcode`O="704F \mathcode`P="7050
\mathcode`Q="7051 \mathcode`R="7052 \mathcode`S="7053 \mathcode`T="7054
\mathcode`U="7055 \mathcode`V="7056 \mathcode`W="7057 \mathcode`X="7058
\mathcode`Y="7059 \mathcode`Z="705A

\mathchardef\nearrow"3225
\mathchardef\searrow"3226
\mathchardef\nwarrow"322D
\mathchardef\swarrow"322E

\newcommand{\hfl}[2]{\smash{\mathop{\hbox to 6mm{\rightarrowfill}}
                     \limits^{\scriptstyle#1}_{\scriptstyle#2}}}

\newcommand{\efl}[2]{\smash{\mathop{\hbox to 6mm{\rightarrowfill}}
                     \limits^{\scriptstyle#1}_{\scriptstyle#2}}}

\newcommand{\sfl}[2]{\llap{$\scriptstyle #1$}\left\downarrow
                     \vbox to 3mm{}\right.\rlap{$\scriptstyle #2$}}

\newcommand{\wfl}[2]{\smash{\mathop{\hbox to 6mm{\leftarrowfill}}
                     \limits^{\scriptstyle#1}_{\scriptstyle#2}}}

\newcommand{\sefl}[2]{\raise2pt\hbox{\llap{$\scriptstyle#1$}}\hspace{-2pt}%
                  \searrow\hspace{-2pt}\raise2pt\hbox{\rlap{$\scriptstyle#2$}}}

\newcommand{\swfl}[2]{\raise2pt\hbox{\llap{$\scriptstyle#1$}}\hspace{-2pt}
                  \swarrow\hspace{-2pt}\raise2pt\hbox{\rlap{$\scriptstyle#2$}}}

\newenvironment{diagram}{\begin{matrix}}{\end{matrix}}

\newcommand{\nc}{\newcommand}
\nc{\on}{\operatorname}
\nc{\onw}{\operatornamewithlimits}

\nc{\reln}{{\mathbb Z}}
\nc{\quot}{{\mathbb Q}}
\nc{\comp}{{\mathbb C}}
\nc{\proj}{{\mathbb P}}
\nc{\real}{{\mathbb R}}
\nc{\natn}{{\mathbb N}}
\nc{\planp}{\proj^2}
\nc{\aff}{{\mathbb A}}
\nc{\tr}{\mathbb}
\nc{\goth}{\mathfrak}
\nc{\scr}{\mathscr}

\nc{\isom}{\stackrel{\sim}{\rightarrow}}
\nc{\lra}{\longrightarrow}
\nc{\ra}{\rightarrow}
\nc{\GL}{\on{GL}}
\nc{\SL}{\on{SL}}
\nc{\U}{\on{U}}
\nc{\SO}{\on{SO}}
\nc{\Spin}{\on{Spin}}
\nc{\Sym}{\on{Sym}}
\nc{\Res}{\onw{Res}}
\nc{\Ext}{\on{Ext}}
\nc{\Tor}{\on{Tor}}
\nc{\Hom}{\on{Hom}}
\nc{\End}{\on{End}}
\nc{\Pic}{\on{Pic}}
\renewcommand{\H}{\on{H}}
\nc{\Todd}{\on{Todd}}
\nc{\res}{\on{Res}}
\nc{\Spec}{\on{Spec}}
\nc{\Gr}{\on{Gr}}
\nc{\Ind}{\on{Ind}}
\nc{\Vac}{\on{Vac}}
\nc{\Tr}{\on{Tr}}
\nc{\tensor}{\otimes}
\nc{\rang}{\on{rang}}
\nc{\ds}{\displaystyle}
\nc{\supp}{\on{supp}}
\nc{\Quot}{\on{Quot}}
\nc{\supps}{\supp_{s}}
\nc{\prof}{\on{prof}}

\nc{\gr}{\on{gr}}
\nc{\id}{\on{id}}
\nc{\HyperExt}{\on{\tr{E}}xt}
\nc{\HyperH}{\on{\tr{H}}}

\nc{\Ker}{\on{Ker}}
\nc{\Coker}{\on{Coker}}
\nc{\Aut}{\on{Aut}}
\nc{\codim}{codim}
\nc{\Cl}{\on{Cl}}
\nc{\Pf}{\on{Pf}}
\nc{\pf}{\on{pf}}
\nc{\ext}{\on{ext}}
\nc{\h}{\on{h}}
\nc{\mult}{\on{mult}}
\nc{\Groth}{\on{Groth}}
\nc{\Stab}{\on{Stab}}
\nc{\droitep}{\proj_{1}}
\nc{\ie}{{\em i.e. }}
\nc{\ul}{\underline}
\nc{\ol}{\overline}

\nc{\inject}{\hookrightarrow}
\nc{\cf}{{\em cf.} }
\nc{\ort}{\perp}
\nc{\R}{\on{R}}
\nc{\tvi}{\vrule height 12pt depth 5pt width 0pt}
\nc{\tv}{\tvi\vrule}
\nc{\osum}{\onw{\oplus}}
\nc{\np}{\clearpage}
\nc{\Vir}{\on{Vir}}
\renewcommand{\isom}{\stackrel{\sim}{\,\ra}}
\renewcommand{\thesubsubsection}{(\thesubsection.\arabic{subsubsection})}
\nc{\no}{\mathop{\raisebox{-2.5pt}{$\stackrel{\textstyle\circ}{\circ}$}}}
\nc{\ensfrac}[2]{\raisebox{+4pt}{$#1$}/\raisebox{-4pt}{$#2$}}
\nc{\motimes}{\dot{\otimes}}
\nc{\moins}{\mathrel{\hbox{\vrule height 3pt depth -2pt width 6pt}}}
\nc{\moinss}{\mathrel{\hbox{\kern1pt\vrule height 2.3pt depth -1.6pt width
4.2pt\kern1pt}}}
\nc{\cqfd}{$\square$}
\nc{\ob}{\on{ob}}
\nc{\car}{\on{char}}
\nc{\Lie}{\on{Lie}}
\nc{\Isom}{\on{Isom}}
\nc{\Mor}{\on{Mor}}
\nc{\extern}{\onw{\boxtimes}}
\nc{\Ad}{\on{Ad}}
\nc{\rank}{\on{rank}}
\nc{\bk}{\backslash}
\nc{\lla}{\longleftarrow}
\nc{\limind}{\mathop{\oalign{lim\cr\hidewidth$\lra$\hidewidth\cr}}}
\nc{\limproj}{\mathop{\oalign{lim\cr\hidewidth$\lla$\hidewidth\cr}}}
\nc{\mono}{\hookrightarrow}
\nc{\omal}{\onw{\otimes}}

\nc{\restriction}[1]{_{|#1}}

\nc{\comment}[1]{}

\renewcommand{\Pf}{\on{Pf}}
\renewcommand{\det}{\on{det}}
\renewcommand{\codim}{\on{codim}}

\nc{\Det}{\on{Det}}
\nc{\Id}{\on{Id}}
\nc{\ev}{\on{ev}}
\nc{\idbb}{{\mathbb I}}
\let\scr\cal

\DeclareMathSymbol{\boxtimes}{\mathbin}{AMSa}{"02}
\nc{\g}{{\goth{g}}}
\nc{\Lg}{L_{\hbox{\normalsize${\goth{g}}$}}}
\nc{\Lgmm}{\Lg^{--}}
\nc{\Lgpp}{\Lg^{++}}
\nc{\Lgh}{\widehat{\Lg}}
\nc{\Lgp}{\Lg^{+}}
\nc{\Lgph}{\widehat{\Lgp}}
\nc{\LgX}{\Lg^{X}}
\nc{\LGp}{L^{+}_{G}}
\nc{\LGm}{L^{-}_{G}}
\nc{\LGmm}{L^{--}_{G}}
\nc{\LG}{L_G}
\nc{\LGh}{\widehat{\LG}}
\nc{\LGX}{\LG^{X}}
\nc{\LGXi}{\LG^{X,i}}
\nc{\LGXone}{\LG^{X,1}}
\nc{\LGXip}{\LG^{X,i+1}}
\nc{\LGXN}{\LGX(N)}
\nc{\LGXiN}{\LGXi(N)}
\nc{\LGN}{\LG(N)}
\nc{\LSL}{LSL_{r}}
\nc{\LSLh}{\widehat{\LSL}}
\nc{\LSLp}{L^{+}_{SL_{r}}}
\nc{\LSLX}{L^{X}_{SL_{r}}}
\nc{\Lsl}{L_{\hbox{\normalsize${\goth{sl}}_{r}$}}}
\nc{\LslV}{L_{\hbox{\normalsize${\goth{sl}}(V)$}}}
\nc{\LslVh}{\widehat{\LslV}}
\nc{\Lslh}{\widehat{\Lsl}}
\nc{\rond}{\circ}
\nc{\G}{{\bf G}}
\nc{\HG}{{\cal{H}}}
\nc{\M}{\on{Bun}_G}
\nc{\MSpin}{\on{Bun}_{\Spin_{r}}}
\nc{\MSO}{\on{Bun}_{SO_{r}}}
\nc{\MGtwo}{\on{Bun}_{G_2}}
\nc{\MSOzero}{\on{Bun}_{SO_{r}}(0)}
\nc{\MSOone}{\on{Bun}_{SO_{r}}(1)}
\nc{\Mproj}{\on{Bun}_G(\proj^{1}_{k})}
\nc{\MM}{\on{Bun}_G^{\prime})}
\nc{\Q}{{\cal{Q}}_G}
\nc{\QSL}{{\cal{Q}}_{SL_{r}}}
\nc{\QN}{{\cal{Q}}_G(N)}
\nc{\Qzero}{{\cal{Q}}_G(0)}
\nc{\QNred}{{\cal{Q}}_G(N)_{\on{red}}}
\nc{\Qp}{{\cal{Q}}_G^{p_{i+1}}}
\nc{\Qpar}{{\cal{Q}}^{par}_G(\ul{p},\ul{P})}
\nc{\Mpar}{\on{Bun}^{par}_G(\ul{p},\ul{P})}
\nc{\ModSO}{M_{SO_{r}}}
\nc{\ModSpin}{M_{\Spin_{r}}}
\nc{\ModGLzero}{M_{GL_{r}}(0)}
\nc{\Modzero}{M_{G}(0)}
\nc{\ModSOreg}{M_{SO_{r}}^{reg}}
\nc{\ModSOone}{M_{SO_{r}}(1)}
\nc{\ModSOzero}{M_{SO_{r}}(0)}
\nc{\Modt}{M_{G}(\tau)}
\nc{\Mod}{M_{G}}
\nc{\Modreg}{M_{G}^{reg}}
\nc{\Modtreg}{M^{reg}_{G}(\tau)}
\nc{\Mt}{\on{Bun}_G(\tau)}
\nc{\Rt}{R_G(\tau)}
\nc{\Mtss}{\on{Bun}^{ss}_{G}(\tau)}
\nc{\Mss}{\on{Bun}^{ss}_{G}}
\nc{\Mreg}{\on{Bun}^{reg}_{G}}
\nc{\MSL}{\on{Bun}_{SL_{r}}}
\nc{\MGL}{\on{Bun}_{GL_{r}}}
\nc{\LGmmm}{L_{-}G_{m}}
\nc{\ma}{{\goth{M}}}
\nc{\Quad}{\on{Quad}}

\nc{\yladress}{ Yves {\sc Laszlo}\\ {\sc  DMI, \'Ecole normale sup\'erieure,
45, rue d'Ulm,  75230 Paris Cedex 05, France}\\ {\it e-mail:}
yves.laszlo@@ens.fr}

\nc{\csadress}{ Christoph {\sc Sorger}\\ {\sc Institut de math\'ematiques
Jussieu (UMR 9994 du CNRS) - Case Postale 7012,  Universit\'e Paris 7, 2, place
Jussieu, F-75251 Paris Cedex 05, France}\\ {\it e-mail:}
sorger@@mathp7.jussieu.fr}

\begin{document}

\title{The line bundles on the moduli of parabolic \\$G$-bundles over curves
and
their sections}
\author{Yves Laszlo and Christoph Sorger}
\date{}
\maketitle

\begin{abstract}

Let $X$ be a smooth, complete and connected curve and $G$ be a simple and
simply
connected algebraic group over $\comp$. We calculate the Picard group of the
moduli stack of quasi-parabolic $G$-bundles and identify the spaces of
sections of its members to the conformal blocs of Tsuchiya, Ueno and Yamada.
We describe the canonical sheaf on these stacks and show that they admit a
unique square root, which we will construct explicitly. Finally we  show how
the
results on the stacks apply to the coarse moduli spaces and recover (and
extend) the Drezet-Narasimhan theorem. We show moreover that the coarse moduli
spaces  of semi-stable
$SO_r$-bundles are not locally factorial for $r\geq 7$.

\end{abstract}

\section{Introduction.}

\subsection{}\label{th:Pic} Fix a simple and simply connected algebraic group
$G$ over $k=\comp$ and a Borel subgroup $B\subset G$. Let $X$ be a smooth,
complete and connected curve over $k$ and
$p_{1},\dots,p_{n}$ be distinct points of $X$, labeled by standard (\ie
containing
$B$) parabolic subgroups $P_{1},\dots,P_{n}$ of $G$ (we allow $n=0$). Let
$\Mpar$ be the moduli stack of quasi-parabolic
$G$-bundles of type $\ul{P}=(P_{1},\dots,P_{n})$
at $\ul{p}=(p_{1},\dots,p_{n})$ and denote by
$X(P_{i})$ the character group of $P_{i}$.

\begin{th*}{} There is
${\scr{L}}\in\Pic(\M)$ such that we have an isomorphism
$$\gamma:\Pic(\Mpar)\isom\reln{\scr{L}}\times\prod_{i=1}^{n}X(P_{i}).$$ If $G$
is of type $A$ or $C$, then ${\scr{L}}$ is the  determinant of cohomology
(\cf\ref{the-det-bundle}). If $G$ is of type $B$, $D$ or $G_{2}$, then
${\scr{L}}$ is the pfaffian of cohomology (\cf\ref{the-pf-bundle}).
\end{th*}

If $G$ is of type $E_{6},E_{7},E_{8},F_{4}$ we believe that we have
${\scr{L}}^{\otimes d(G)}={\scr{D}}_{\rho(G)}$ where respectively
$d(G)=6,12,60,6$ and
${\scr{D}}_{\rho(G)}$ is the determinant of cohomology
(\cf\ref{the-det-bundle}) associated to the fundamental representation
$\rho(G)=\varpi_{6},\varpi_{7},\varpi_{8},\varpi_{4}$
(\cf the discussion in \ref{pb}).

\subsection{} Suppose that the points $p_{1},\dots,p_{n}$ are instead labeled
by finite dimensional simple representations $\lambda_{1},\dots,\lambda_{n}$ of
$G$ and that an additional integer $\ell$, the {\em level}, is fixed.  The
choice of a representation
$\lambda$ of $G$ is equivalent to the choice of a standard parabolic subgroup
$P\subset G$ and a character $\chi\in X(P)$.  Therefore, the labeling of the
points
$p_{1},\dots,p_{n}$ by the representations
$\lambda_{1},\dots,\lambda_{n}$ defines the type $\ul{P}$ of a quasi-parabolic
$G$-bundle,  that is the stack $\Mpar$ {\em and}, by the above theorem, a line
bundle
${\scr{L}}(\ell,\ul\chi)$ over $\Mpar$.  The global sections of
${\scr{L}}(\ell,\ul\chi)$ give a vector space, the space of {\em generalized
parabolic $G$-theta-functions of level $\ell$},
 which is canonically associated to
$(X,\ul{p},\ul\lambda)$. In mathematical physics, the rational conformal field
theory of Tsuchiya, Ueno and Yamada \cite{TUY} associates also to
$(X,\ul{p},\ul\lambda,\ell)$ a vector space: the space of {\em conformal
blocks}
$V_{X}(\ul{p},\ul\lambda,\ell)$ (\cf \cite{So3} for an overview).
\begin{th*}{}
Suppose that $G$ is  classical or $G_2$. Then there is a canonical
isomorphism
\begin{equation}\label{Verlinde} H^{0}(\Mpar,{\scr{L}}(\ell,\ul\chi))\isom
V_{X}(\ul{p},\ul\lambda,\ell).
\end{equation}
In particular, $\dim H^{0}(\Mpar,{\scr{L}}(\ell,\ul\chi))$ is given by the
Verlinde dimension formula.
\end{th*}

For $n=0$, this has been proved independently by Beauville and
the first author
\cite{BL1} for $G=SL_{r}$ and by Faltings \cite{F} and Kumar,
Narasimhan and Ramanathan \cite{KNR} for arbitrary
simple and simply connected $G$. For arbitrary $n$ and $G=SL_r$ this has been
proved by Pauly \cite{P} and can be proved for arbitrary simple and simply
connected $G$ using (\ref{th:Uniformization}) (and therefore \cite{DS}) and
(\ref{LGX-is-integral}) below
 following the lines of
\cite{BL1} and \cite{P}. This is the subject of Section \ref{Identification}.

\subsection{}\label{th:Uniformization}  The above results are proved via the
{\em uniformization} theorem: restrict for simplicity of the introduction to
$n=0$. Suppose $p\in X$ and denote $X^{*}=X\moins p$. Define
$D=\Spec(\hat{\cal{O}}_{p})$, where $\hat{\cal{O}}_{p}$ is the formal
completion of the local ring ${\cal{O}}_{p}$ at $p$ and $D^{*}=\Spec(K_{p})$
where $K_{p}$ is the quotient field of $\hat{\cal{O}}_{p}$. Let $\LG$ (resp.
$\LGp$, resp. $\LGX$) be the group of algebraic morphisms from $D^{*}$ (resp.
$D$, resp. $X^{*}$) to $G$.

\begin{th*}{} The algebraic stack $\M$  is canonically isomorphic to the double
quotient stack $\LGX\backslash\LG/\LGp$. Moreover, the projection map
$$\pi:\Q:=\LG/\LGp\ra \M$$ is locally trivial for the
\'etale topology.
\end{th*}

This is proved in \cite{BL1} for $G=SL_{r}$. The extension to arbitrary $G$ has
been made possible by Drinfeld and Simpson
\cite{DS} in response to a question by the first author. They prove  that if
$S$
is a
$k$-scheme and
$E$ a
$G$-bundle over
$X\times S$ then, locally for the \'etale  topology on $S$, the restriction of
$E$ to
$X^{*}\times S$ is trivial, which is essential for the proof.  The above
theorem is valid more generally for semi-simple $G$. In characteristic $p$, one
has to replace ``\'etale'' by ``fppf'' if $p$ divides the order of
$\pi_{1}(G(\comp))$.

\subsection{}\label{pb} Consider the pullback  morphism
$$\pi^{*}:\Pic(\M)\efl{}{}\Pic(\Q).$$

The Picard group of $\Q$ is known (\cite{Ma},\cite{KNR}) to be canonically
isomorphic to $\reln$, which reduces proving Theorem \ref{th:Pic} to proving
that
$\pi^{*}$ is an isomorphism. We will show that the injectivity of $\alpha$ will
follow from  the fact that $\LGX$ has no characters which in turn will follow
from the fact that $\LGX$ is reduced and connected.  Moreover, the surjectivity
of $\alpha$ would follow from the simple connectedness of $\LGX$.  Both
topological properties, connectedness and simple connectedness of $\LGX$ are
affirmed in \cite{KNR} and we believe them to be true.
Whereas we will prove the connectedness of $\LGX$,
following an idea of V. Drinfeld, we do not see how to prove the simple
connectedness of $\LGX$.
The injectivity is enough to prove the first part of Theorem \ref{th:Pic},
but to identify the generator ${\cal{L}}$ we should prove the surjectivity of
$\alpha$. For classical $G$ and
$G_{2}$ we do this by constructing a line bundle on $\M$ pulling back to the
generator of $\Pic(\Q)$.

\subsection{}\label{sr} The pfaffian construction (\ref{the-pf-bundle}) may be
used to prove the following, valid over $k$ algebraically closed of
characteristic $\not=2$.
\begin{prop*}{} Suppose $G$ is semi-simple. Then, for every
theta-characteristic
$\kappa$ on $X$, there is a canonical square-root ${\cal{P}}_{\kappa}$ of the
dualizing sheaf
$\omega_{_{\M}}$ of $\M$.
\end{prop*}

\subsection{}\label{th:Pic(Mmod)} The last section will be devoted to
Ramanathans moduli spaces $\Mod$ of semi-stable
$G$-bundles.  We will show how some of the results for the stack $\M$ will be
true also for the moduli spaces $\Mod$. In particular we will recover (and
extend) the Drezet-Narasimhan  theorem.

\begin{th*}{} There is a canonical isomorphism
$\Pic(\Mod)\isom\reln L.$ If $G$ is of type $A$ or $C$ then $L$ is the
determinant bundle and moreover, in this case $\Mod$ is locally factorial. If
$G$ is of type
$B_r$, $D_r$, $r\geq 4$ or
$G_2$ then $L$ or $L^{\otimes 2}$ is the determinant bundle.
\end{th*}

This theorem has also been proved, independently and with a different
method, by Kumar and Narasimhan \cite{KN}.

The question whether $\Mod$ is locally factorial for $G$ of type other than
the simply connected groups of type $A,C$ is the subject of a forthcoming
paper. We show there for example that $\Pic(\ModSpin)$ is generated
by the determinant of cohomology  and in particular that
$\ModSpin$ is not locally factorial for $r\geq 7$ by
``lifting'' to $\Spin_r$ the proof we give here
(\ref{not-locally-factorial}) for the analogous statement for
$\ModSO$.

\bigskip

{\em We would like to thank A. Beauville and C. Simpson for useful
discussions and V. Drinfeld for his suggestion in (\ref{LGX-is-integral}) 
and for pointing out an inaccuracy in an earlier version of this paper.}

\section{Some Lie theory.}\label{Lie-theory}

Throughout this section $k$ will be an algebraically closed field of
characteristic zero.

\subsection{The general set up.} Let ${\goth{g}}$ be a simple finite
dimensional Lie algebra over $k$. We fix a Cartan subalgebra
${\goth{h}}\subset{\goth{g}}$ and denote by $\Delta$ the associated root
system. We have the root decomposition
$\displaystyle{\goth{g}}={\goth{h}}\oplus\osum_{\alpha\in
\Delta}{\goth{g}}_{\alpha}.$   The Lie subalgebra
${\goth{g}}_{-\alpha}\oplus[{\goth{g}}_{\alpha},{\goth{g}}_{-\alpha}]\oplus
{\goth{g}}_{\alpha}$, isomorphic as a Lie algebra to ${\goth{sl}}_{2}$, will be
denoted by ${\goth{sl}}_{2}(\alpha)$. Moreover we choose a basis
$\Pi=\{\alpha_{1},\dots,\alpha_{r}\}$ of $\Delta$ and we denote by
$\Delta_{+}$ the set of positive roots (with respect to $\Pi$). Put
$\ds{\goth{b}}={\goth{h}}\oplus(\osum_{\alpha\in
\Delta_{+}}{\goth{g}}_{\alpha})$. For each
$\alpha\in \Delta_{+}$, we denote by
$H_{\alpha}$ the coroot of
$\alpha$, \ie the unique element of
$[{\goth{g}}_{\alpha},{\goth{g}}_{-\alpha}]$ such that
$\alpha(H_{\alpha})=2$, and we denote by
$X_{\alpha}\in{\goth{g}}_{\alpha}$ and
$X_{-\alpha}\in{\goth{g}}_{-\alpha}$ elements such that
$$[H_{\alpha},X_{\alpha}]=2X_{\alpha},\hspace{1cm}
[H_{\alpha},X_{-\alpha}]=-2X_{-\alpha},\hspace{1cm}
[X_{\alpha},X_{-\alpha}]=H_{\alpha}.$$ When $\alpha$ is one of the simple roots
$\alpha_{i}$, we write
$H_{i},X_{i},Y_{i}$ instead of $H_{\alpha_{i}},X_{a_{i}},Y_{\alpha_{i}}$. Let
$(\varpi_{i})$ be the basis of ${\goth{h}}^{*}$ dual to the basis
$(H_{i})$. Denote by $P$ the weight lattice and by
$P_{+}\subset P$ the set of dominant weights. Given a dominant weight
$\lambda$, denote
$L_{\lambda}$ the associated simple ${\goth{g}}$-module with highest weight
$\lambda$ and $v_{\lambda}$ its highest weight vector. Finally
$(\,,\,)$ will be the Cartan-Killing form normalized such that for the highest
root
$\theta$ we have $(\theta,\theta)=2$.

\subsection{Loop algebras.} Let $\Lg={\goth{g}}\otimes_{k}k((z))$ be the {\em
loop algebra} of
$\g$ and $\Lgp={\goth{g}}\otimes_{k}k[[z]]$ its subalgebra of {\em positive}
loops.
 There is a natural $2$-cocycle
$$\begin{diagram}
\psi_{\g}:&\Lg&\times&\Lg&\lra&k\\ &(X\otimes f&,&Y\otimes
g)&\mapsto&(X,Y)\Res(gdf)\\
\end{diagram}
$$ defining a central extension $\Lgh$ of $\g$:
$$0\lra k\lra\Lgh\lra\Lg\lra 0.$$ Every other cocycle is a scalar multiple of
$\psi_{\g}$ and the above central extension is universal. Let $\Lgph$ be the
extension of $\Lgp$ obtained by restricting the above extension to $\Lgp$. As
the cocycle is trivial over
$\Lgp$ this extension splits.

Let $\ell$ be a positive integer. A representation of $\Lgh$ is {\em of level
$\ell$} if the center $c$ acts by multiplication by $\ell$. Such a
representation is called integrable if $X\tensor f$ acts locally nilpotent for
all $X\tensor f\in{\goth{g}}_{\alpha}\otimes_{k}k((z))$. The theory of affine
Lie algebras
\cite{Kac} affirms that the irreducible integrable representations of level
$\ell$ of
$\Lgh$ are classified (up to isomorphism) by the weights
$P_{\ell}=\{\lambda\in P_{+}/(\lambda,\theta)\leq\ell\}.$ We denote by
${\cal{H}}_{\lambda,\ell}$ the irreducible integrable representation of level
$\ell$ and highest weight $\lambda\in P_{\ell}$. If $\lambda=0$, the
corresponding representation, which we denote simply by
${\cal{H}}_{\ell}$, is called the {\em basic representation of level $\ell$}.

\subsection{The Dynkin index.}\label{section-Dynkin-index} Let
$\rho:\g\ra{\goth{sl}}(V)$ be a representation of $\g$. Then $\rho$ induces a
morphism of Lie algebras
$\Lg\ra\LslV$. Pull back the universal central extension to
$\Lg$:
$$
\begin{diagram} 0&\lra&k&\lra&\widetilde{\Lg}&\lra&\Lg&\lra&0\\
 &&\parallel&&\sfl{}{}&&\sfl{}{}\\
0&\lra&k&\lra&\LslVh&\lra&\LslV&\lra&0\\
\end{diagram}
$$ The cocycle of the central extension $\widetilde{\Lg}$ is of the type
$d_{\rho}\psi_{\g}$. Define the {\em Dynkin index} of the representation
$\rho$ of $\g$ by the number $d_{\rho}$.

\begin{lem} Let $V=\sum_{\lambda}n_{\lambda}e^{\lambda}$ be the formal
character of
$V$. Then we have
$$d_{\rho}=\frac{1}{2}\sum_{\lambda}n_{\lambda}\lambda(H_{\theta})^{2}$$
\end{lem}

\begin{proof} By definition of the cocycle, we have
$d_{\rho}=\Tr(\rho(X_{\theta})\rho(X_{-\theta}))$. Decompose the
${\goth{sl}}_{2}(\theta)$-module,
$V$ as $\osum V^{(d_i)}$, where $V^{(d_i)}$ is the standard irreducible
${\goth{sl}}_{2}$-module with highest weight $d_i$. We may realize $V^{(d_i)}$
as the vector space of homogeneous polynomials in 2 variables
$x$ and $y$ of degree $d_i$. Then $X_\theta$ acts as $x\partial/\partial y$,
and
$X_{-\theta}$ as $y\partial/\partial x$. Using the basis $x^ly^{d_i-l},
l=0,\ldots, d_i$ of $V(d_i)$, we see
$$d_{\rho}=\sum_i\sum_{k=0}^{d_i}k(d_i+1-k).$$ The formal character of the
${\goth{sl}}_{2}(\theta)$-module $V^{(d)}$ is
$\sum_{k=0}^{d}e^{d\rho_{\theta}-k\alpha_{\theta}}$ where $\alpha_{\theta}$ is
the positive root of ${\goth{sl}}_{2}(\theta)$ and
$\rho_{\theta}=\frac{1}{2}\alpha_{\theta}$. Therefore we are reduced to prove
the equality
$$\sum_{k=0}^{d}k(d+1-k)=
\frac{1}{2}\sum_{k=0}^{d}[d\rho_{\theta}-k\alpha_{\theta})(H_{\theta})]^{2}
=\frac{1}{2}\sum_{k=0}^{d}[d-2k]^{2}
$$ which is easy.
\end{proof}

\begin{rem} The Dynkin index of a representation has been introduced to the
theory of
$G$-bundles over a curve by Faltings \cite{F} and Kumar, Narasimhan, Ramanathan
\cite{KNR}.
\end{rem}

We are interested here in the {\em minimal Dynkin index}
$d_{{\goth{g}}}$ defined to be as $\min d_{\rho}$ where
$\rho$ runs over all representations $\rho:\g\ra{\goth{sl}}(V)$.

\begin{prop}\label{minimal-Dynkin-index} The minimal Dynkin index
$d_{{\goth{g}}}$ is as follows
$$\begin{array}{c|c|c|c|c|c|c|c|c|c}
\text{Type of }{\goth{g}}&A_{r}&B_{r},r\geq 3&C_{r}&D_{r},\geq
4&E_{6}&E_{7}&E_{8}&F_{4}&G_{2}\\\hline
d_{{\goth{g}}}&1&2&1&2&6&12&60&6&2\\\hline
\lambda\text{ s.t. } d_{{\goth{g}}}=d_{\rho(\lambda)}&
\varpi_{1}&\varpi_{1}&\varpi_{1}&\varpi_{1}&\varpi_{6}&\varpi_{7}&\varpi_{8}
&\varpi_{4}&\varpi_{1}\\
\end{array}
$$ Moreover, for any representation $\rho:\g\ra{\goth{sl}}(V)$, we have
$d_{\rho}=0\bmod d_{{\goth{g}}}$.
\end{prop}

{\it Proof.} It is enough to calculate the Dynkin index for the fundamental
weights (note that $d_{V\otimes W}=r_{W}d_{V}+r_{V}d_{W}$ if $V$ and $W$ are
two ${\goth{g}}$-modules of rank
$r_{V}$ and $r_{W}$), which can be done explicitly \cite{D}. We give the values
here for $E_8$, as not all of them in (\cite{D}, Table 5) are
correct:
$$\begin{array}{c|c|c|c|c|c|c|c|c|} &d_{\varpi_{1}}&d_{\varpi_{2}}&
d_{\varpi_{3}}&d_{\varpi_{4}}&d_{\varpi_{5}}&d_{\varpi_{6}}&
d_{\varpi_{7}}&d_{\varpi_{8}}\\\hline
E_{8}&1500&85500&5292000&8345660400&141605100&1778400&14700&60\\\hline
\end{array}
$$

\section{The stack $\M$.}\label{the-stack-M}

Throughout this section $k$ will be an algebraically closed field, $G$ will be
semi-simple algebraic group over $k$.

\subsection{} Let $X$ be a scheme over $k$. By a {\em principal $G$-bundle}
over $X$ (or just $G$-bundle for short), we understand a scheme $E\ra
X$ equipped with a right action of $G$ such that, locally in the flat topology,
$E$ is trivial, \ie isomorphic to $G\times X$ as an $G$-homogeneous space. In
particular, $E$ is affine, flat and smooth over $X$. Moreover,  the above
conditions imply that $E$ is even locally trivial for the \'etale topology.

If $F$ is a quasi-projective scheme on which $G$ acts on the left and $E$ is a
$G$-bundle, we can form $E(F)=E\times^{G}F$ the {\em associated bundle with
fiber
$F$}. It is the quotient of
$E\times F$ under the action of $G$ defined by $g.(e,f)=(e.g,g^{-1}f)$. If $H$
is a subgroup of $G$, the associated $G/H$-bundle $E(G/H)$ will be denoted
simply by $E/H$.

Let $\rho:G\ra G^{\prime}$ be a morphism of algebraic groups. Then, as $G$ acts
on
$G^{\prime}$ via $\rho$, we can form the {\em extension of the structure group}
of a
$G$-bundle $E$, that is the $G^{\prime}$-bundle $E(G^{\prime})$. Conversely, if
$F$ is a $G^{\prime}$-bundle, a {\em reduction of structure group}
$F_{G}$ is a $G$-bundle $E$ together with an isomorphism
$F_{G}(G^{\prime})\isom F$. If $\rho$ is a closed immersion, such reductions
are in one to one correspondence with section of the associated bundle $F/G$.

\subsection{} Let us collect some well known generalities on stacks for further
reference. Let ${\tr A}\text{ff/$k$}$ be the flat affine site over $k$, that is
the category of $k$-algebras equipped with the {\em fppf} topology. By {\em
$k$-space} (resp. {\em $k$-group}) we understand a sheaf of sets (resp. groups)
over
${\tr A}\text{ff/$k$}$. Any $k$-scheme can (and will) be considered as a
$k$-space.

We will view {\em $k$-stacks} from the pseudo-functorial point of view, \ie a
$k$-stack
${\goth{X}}$ will associate to every $k$-algebra
$R$ a groupoid
${\goth{X}}(R)$ and to every morphism of $k$-algebras
$u:R\ra R^{\prime}$ a functor
$u^{*}:{\goth{X}}(R^{\prime})\ra{\goth{X}}(R)$ together with isomorphisms of
functors
$(u\circ v)^{*}\simeq v^{*}\circ u^{*}$ satisfying the usual cocycle condition.
The required topological properties are that for every
$x,y\in\ob{\goth{X}}(R)$ the presheaf $\ul\Isom(x,y)$ is a sheaf and that all
descent data are effective (\cite{LMB}, 2.1). Any $k$-space $X$ may be seen as
a
$k$-stack, by considering a set as a groupoid (with the identity as the only
morphism). Conversely, any $k$-stack ${\goth{X}}$ such that ${\goth{X}}(R)$ is
a
discrete groupoid (\ie has only the identity as automorphisms) for all
$k$-algebras
$R$, is a $k$-space.

A morphism $F:{\goth{X}}\ra{\goth{Y}}$ will associate, for every $k$-algebra
$R$, a functor ${\goth{X}}(R)\ra{\goth{Y}}(R)$ satisfying the obvious
compatibility conditions. Let $S=\Spec(R)$ and consider a morphism
$\eta:S\ra{\goth{Y}}$, that is an object $\eta$ of ${\goth{Y}}(S)$. The fiber
${\goth{X}}_{\eta}$ is a stack over
$S$. The morphism $F$ is {\em representable} if
${\goth{X}}_{\eta}$ is representable as a scheme for all $S=\Spec(R)$. A stack
${\goth{X}}$ is {\em algebraic} if the diagonal morphism
${\goth{X}}\ra{\goth{X}}\times{\goth{X}}$ is representable, separated and
quasi-compact and if there is a scheme $X$ and a representable, smooth,
surjective morphism of stacks $P:X\ra{\goth{X}}$.

Suppose $X$ is a $k$-space and that the $k$-group $\Gamma$ acts on $X$. Then
the
quotient stack $[X/\Gamma]$ is defined as follows. Let $R$ be a $k$-algebra.
The
objects of $[X/\Gamma](R)$ are pairs $(E,\alpha)$ where $E$ is a
$\Gamma$-bundle over $\Spec(R)$ and $\alpha:E\ra X$ is
$G$-equivariant, the arrows are defined in the obvious way and so are the
functors
$[X/\Gamma](R^{\prime})\ra [X/\Gamma](R)$.

\subsection{} Let $X$ be a smooth, complete and connected curve of genus $g$
over
$k$. We denote by
$\M$ the  stack of principal $G$-bundles over $X$ which is defined as follows.
For any $k$-algebra $R$ denote $X_R$ the scheme $X\times_{k}\Spec(R)$. Then
objects of
$\M(R)$ are $G$-bundles over $X_R$, morphisms of
$\M(R)$ are isomorphisms of $G$-bundles.

The following proposition is well known.

\begin{prop} The stack $\M$ is algebraic and smooth. Moreover we have
$\dim\M=(g-1)\dim G.$
\end{prop}

\subsection{} Choose a closed point $p$ on $X$ and set $X^*=X\moins p$.  Let
${\cal O}$ be the completion of the local ring of $X$ at $p$, and $K$ its field
of fractions. Set $D=\Spec({\cal{O}})$ and $D^{*}=\Spec(K)$. We choose a local
coordinate $z$ at $p$ and identify ${\cal O}$ with $k[[z]]$ and $K$ with
$k((z))$. Let $R$ be a
$k$-algebra. Define $X_R = X\times_k \Spec(R)$,
$X_R^* = X^*\times_k \Spec(R)$, $D_R = \Spec\bigl(R[[z]]\bigr)$ and
$D_R^* = \Spec\bigl(R((z))\bigr)$. Then we have the cartesian diagram

$$\begin{diagram} D_R^*&\efl{}{}&D_{R}\\
\sfl{}{}&&\sfl{}{}\\ X_{R}^{*}&\efl{}{}&X_{R}\\
\end{diagram}$$

We denote by $A_{X_R}$ the $k$-algebra
$\Gamma(X_{R}^*\, ,{\cal O}_{X^{*}_R})$.

\subsection{Loop groups.}\label{Loop-groups} The category of $k$-spaces is
closed under direct limits. A $k$-space ($k$-group) will be called an {\em
ind-scheme} (resp. {\em ind-group}) if it is direct limit of a directed system
of schemes. Remark that an ind-group is not necessarily an inductive limit of
algebraic groups.

We denote by $\LG$ the {\em loop group} of $G$ that is the $k$-group defined
$R\mapsto G\bigl(R((z))\bigr)$, where $R$ is any $k$-algebra. The group of {\em
positive loops}, that is the $k$-group
$R\mapsto G\bigl(R[[z]]\bigr)$ will be denoted by by $\LGp$ and the group of
{\em negative loops}, that is the $k$-group
$R\mapsto G\bigl(R[z^{-1}]\bigr)$ will be denoted by $\LGm$. The group of {\em
loops coming from $X^{*}$}, \ie the $k$-group defined by
$R\mapsto G(A_{X_{R}})$, will be denoted by $\LGX$. Finally, we will use also
the
$k$-group $\LGmm$ defined by
$R\mapsto G\bigl(z^{-1}R[z^{-1}]\bigr)$.

Choose a faithful representation
$G\subset SL_{r}$. For $N\ge 0$, we denote by
$\LGN(R)$ the set of matrices $A(z)$ in $G\bigl(R((z))\bigr)\subset
SL_{r}\bigl(R((z))\bigr)$ such that for both $A(z)$ and $A(z)^{-1}$, the
coefficients have  a pole of order
$\le N$. This defines a subfunctor  $\LGN$ of  $\LG$ which is obviously
representable by an (infinite dimensional) affine $k$-scheme.

\begin{prop}\label{ind-groups} The $k$-group  $\LGp$ is an affine group scheme.
The
$k$-group
$\LG$ is an ind-group, direct limit of the sequence of the schemes
$(\LGN)_{N\ge 0}$. Moreover, this ind-structure does not depend on the
embedding
$G\subset SL_{r}$.
\end{prop}

The $k$-group $\LGX$ has the structure of an ind-group induced by the one of
$\LG$. The quotient $k$-space $\Q:=\LG/\LGp$ has equally the structure of an
ind-scheme: define $\QN=\LGN/\LGp$ (note that $\LGN$
is stable under right multiplication by $\Qzero=\LGp$).

\subsection{} Consider triples $(E,\rho,\sigma)$ where $E$ is a $G$-bundle on
$X_R$, $\rho : G\times X_R^* \ra E_{|X_R^*}$ a trivialization of
$E$ over $X_R^*$ and $\sigma : G\times D_R\ra E_{|D_R}$ a trivialization of
$E$ over $D_R$.  We let $T(R)$ be the set of isomorphism
 classes of triples $(E,\rho,\sigma)$.

\begin{prop}\label{triples} The ind-group $\LG$ represents the functor $T$.
\end{prop}

\begin{proof} Let $(E,\rho,\sigma)$ be an element of $T(R)$. Pulling back the
trivializations
$\rho$ and $\sigma$ to $D^*_R$ provides two trivializations
$\rho^*$ and $\sigma^*$ of the pull back of $E$ over ${D^*_R}$: these
trivializations differ by an element $\gamma = \rho^{*-1}\circ\sigma^*$ of $
G\bigl(R((z))\bigr)$.

Conversely, let us start from an element $\gamma$ of $G\bigl(R((z))\bigr)$.
This
element defines an isomorphism of the pullbacks over $D^*_R$ of the trivial
$G$-bundle ${\cal{F}}$ over $X^*_R$ and the trivial $G$-bundle ${\cal{G}}$ over
$D_R$.  These two torsors glue together to a
$G$-bundle $E$ in a functorial way by \cite{BL2} (in fact \cite{BL2} is written
for $SL_r$ but the extension to arbitrary $G$ is straightforward).

These constructions are inverse to each
other by construction.
\end{proof}

\subsection{} Consider the functor $D_G$ which associates to a
$k$-algebra $R$ the set $D_G(R)$ of isomorphism classes of pairs $(E,\rho)$,
where
$E$ is a $G$-bundle over $X_R$ and $\rho$ a trivialization of $E$ over $X^*_R$.

\begin{prop}\label{pairs} The ind-scheme $\Q$ represents the functor $D_{G}$.
\end{prop}

\begin{proof} Let
$R$ be a $k$-algebra and $q$  an element of $\Q(R)$. By definition there exists
a faithfully flat homomorphism $R\rightarrow R'$ and an element $\gamma$ of
$G\bigl(R'((z))\bigr)$ such that the image of $q$ in $\Q(R')$ is the class of
$\gamma$.  To $\gamma$ corresponds by Proposition \ref{triples} a triple
$(E',\rho',\sigma')$ over $X_{R'}$. Let $R''=R'\otimes_R R'$, and let
$(E''_1,\rho''_1)$, $(E''_2,\rho''_2)$ denote the pull-backs of $(E',\rho')$ by
the two projections of $X_{R''}$ onto $X_{R'}$. Since the two images of
$\gamma$ in $G\bigl(R''((z))\bigr)$ differ by an element of
$G\bigl(R''[[z]]\bigr)$, these pairs  are isomorphic. So the isomorphism
$\rho^{\prime\prime}_2 \rho_1^{\prime\prime -1}$ over $X_{R''}^*$ extends to an
isomorphism
$u:E''_1\rightarrow E''_2$ over $X_{R''}$, satisfying the usual cocycle
condition (it is enough to check this over
$X^*$, where it is obvious). Therefore $(E',\rho')$ descends to a pair
$(E,\rho)$ on $X_R$ as in the statement of the proposition.

Conversely, given a pair $(E,\rho)$ as above over $X_R$, we can find a
faithfully flat homomorphism $R\rightarrow R'$ and a trivialization $\sigma'$
of the pull back of $E$ over $D_{R'}$ (in fact, we can take
$R'$ to be the product of henselization of each localized ring $R_x,\
x\in\Spec(R)$). By Proposition \ref{triples} we get an element $\gamma'$ of
$G\bigl(R'((z))\bigr)$ such that the two images of
$\gamma'$ in $G\bigl(R''((z))\bigr)$ (with $R''=R'\otimes_R R'$) differ by an
element of $G\bigl(R''[[z]]\bigr)$; this gives an element of $\Q(R)$. The two
constructions are clearly inverse one of each other.
\end{proof}

We will make use of the following theorem
\begin{th}\label{Drinfeld-Simpson} (Drinfeld-Simpson \cite{DS}) Let $E$ be a
$G$-bundle over $X_R$. Then the restriction of $E$ to $X^*_R$ is trivial {\it
fppf} locally over $\Spec(R)$. If $char(k)$ does not divide the order of
$\pi_{1}(G(\comp))$, then this is even true {\it \'etale} locally.
\end{th}

\subsection{Proof of Theorem \ref{th:Uniformization}.} The universal $G$-bundle
$E$ over $X\times\Q$ (Proposition \ref{pairs}), gives rise to a map
$\pi:\Q\rightarrow \M$. This map is
$\LGX$-invariant, hence induces a morphism of stacks
$\overline{\pi}:\LGX\bk\Q\rightarrow \M$.

On the other hand we can define a map $\M\ra\LGX\bk\Q$ as follows. Let $R$ be a
$k$-algebra, $E$ a $G$-bundle over $X_R$. For any
$R$-algebra $R^{\prime}$, let $T(R^{\prime})$ be the set of trivializations
$\rho$ of
$E_{R^{\prime}}$ over $X^*_{R^{\prime}}$.  This defines a $R$-space $T$ on
which the group $\LGX$ acts.  By Theorem \ref{Drinfeld-Simpson}, it is a torsor
under
$\LGX$. To any element of
$T(R^{\prime})$ corresponds a pair $(E_{R^{\prime}},\rho)$, hence by
Proposition \ref{pairs} an  element of
$\Q(R^{\prime})$. In this way we associate functorially to an object $E$ of
$\M(R)$ a
$\LGX$-equivariant map
$\alpha:T\rightarrow \Q$. This defines a morphism of stacks $\M\rightarrow
\LGX\bk\Q$ which is the inverse of $\overline{\pi}$.

The second assertion means that for any scheme $S$ over $k$ (resp. over $k$
such that
$char(k)$ does not divide the order of $\pi_{1}(G(\comp))$) and any morphism
$f:T\rightarrow \M$, the pull back to $S$  of the fibration
$\pi$ is {\it fppf} (resp. {\it\'etale}) locally trivial, i.e. admits local
sections (for the {\it fppf} (resp. {\it \'etale}) topology). Now $f$
corresponds to a
$G$-bundle $E$ over
$X\times S$. Let $s\in S$. Again by Theorem \ref{Drinfeld-Simpson}, we can find
an {\it fppf} (resp. {\it \'etale}) neighborhood $U$ of $s$ in  $S$ and a
trivialization
$\rho$ of
$E_{|X^*\times U}$. The pair $(E,\rho)$ defines  a morphism $g:U\rightarrow \Q$
(Proposition \ref{pairs}) such that $\pi\rond g=f$, that is a section over $U$
of the pull back of the fibration $\pi$. \cqfd

\section{The infinite grassmannian $\Q$}

Let $G$ be semi-simple and $k$ be an algebraically closed field of
characteristic
$0$ in (\ref{ind-structures-are-the-same}).

\subsection{}\label{mu-is-immersion}
We will use the following two facts:

$(a)$ We may write $\Q$ as direct limit of projective
finite-dimensional $k$-schemes.

$(b)$ The multiplication map
 $\mu:\LGmm\times \LGp  \longrightarrow \LG$ is an open immersion.

\noindent For the first statement, remark that it is enough to consider the
$k$-space
${\cal{Q}}({\goth{g}})$ parameterizing isomorphism classes of sheaves of Lie
algebra which are locally of the form  $S\times{\rm Lie(G)}$ together with  a
trivialization over $\droitep^*$  (note that, looking at the adjoint group,
$\Q$ may be seen as a connected component of ${\cal{Q}}({\goth{g}})$) then
argue as in \cite{BL1}.  For the second statement, the argument of (\cite{BL1}
Proposition 1.11) generalizes to arbitrary $G$, once we know the following.

{\em Claim.} Suppose $Y$ is a proper $S$-scheme and that the structural
morphism has a section $\sigma: S\ra Y$. Suppose moreover that $G\bk H$ is a
reductive subgroup of a reductive group $H$. Then, for any
$G$-bundle $P$ trivial along $\sigma$ the following is true: if the associated
$H$-bundle $P(H)$ is trivial, the $G$-bundle $P$ is so.

Indeed, by assumption, there exists a section $\tau: Y\ra P(H)$. The quotient
$G\bk H$ is affine. Therefore, the composite morphism from $Y$ to $G\bk H$
(which is the composition of
$\tau$ and of the canonical projection $P(H)\ra G\bk H$ factors as
$Y\ra S\hfl{p}{} G\setminus H.$ After an eventual translation of $\tau$ by an
element of $H(S)$, we can assume that the restriction $\sigma^*(\tau)$ of
$\tau$
along $\sigma$ is induced by the trivialization of $P$. Therefore, the morphism
$p$ is the constant morphism with constant value $G\in G\bk H$. In other words,
locally for the \'etale topology on $Y$, the section $\tau$ can be written
$$\tau(y)=\bigl(p(y),h(y)\bigr) {\rm mod}\ G\ {\rm where\ }y\in Y\ {\rm and}\
p(y)\in P, h(y)\in G.$$ The expression $p(y)h(y)^{-1}$ is well defined and
defines a section of $P$.

\subsection{}
\label{ind-structures-are-the-same}
The quotient $\LG/\LGp$ has also been studied by Kumar and Mathieu. But the
structure of ind-variety they put on the quotient is, a priori, not the same as
the functorial one of section \ref{the-stack-M}. As we will use their results,
we have to identify them.

\begin{prop*}{} The ind-structure of $\Q$ defined in section \ref{the-stack-M}
coincides with the ind-structure of Kumar and Mathieu.
\end{prop*}

\begin{proof} Recall that an ind-scheme is called
{\em reduced} (resp. {\em irreducible, integral}) if it is a direct limit of
an
increasing sequence of reduced (resp. irreducible, integral) schemes. By
Lemma 6.3 of \cite{BL1} an ind-scheme is integral if and only if it is
irreducible and reduced.  According to Faltings
\cite{F} (see \cite{BL1} for the case
$SL_{r}$), the ind-group $\LGm$ is integral.  This may
be seen by looking at
$(\LGm)_{red}$ and using Shavarevich's theorem that a closed immersion of
irreducible ind-affine groups which is an isomorphism on Lie algebras, is an
isomorphism
\cite{Sh}. Note that irreducibility is due to the fact that any element can be
deformed to a constant in $G$; that $\Lie(\LGm)\rightarrow
\Lie(\LGm)_{red}$ is an isomorphism can be seen by using the fact that
$\Lie(G)$ is generated by nilpotent elements.

It follows for semisimple $G$ that $\LGmm$, which is a semidirect product of
$G$ and
$\LGm$, is integral, and furthermore if $G$ is simply connected then
$\Q$ is integral.  Indeed by
(\ref{mu-is-immersion} b) $\LGmm\ra\Q$ is an open immersion hence it is enough
to show that $\Q$ is irreducible. Using that connected ind-groups are
irreducible (\cite{Sh}, Proposition 3) and the quotient morphism $\LG\ra\Q$
we reduce to prove the connectedness of $\LG$ which is well known (and follows
for example from uniformization for $\droitep$ and the
corresponding statement for $\Mproj$ , \cf \cite{DS}).

We are ready to deduce the identification of our ind-structure on $\Q$ with the
one used by Kumar or Mathieu in their generalized Borel-Weil theory. Both Kumar
and Mathieu define the structure of ind-variety on
$\LG/\LGp$  using representation theory of Kac-Moody algebras; for instance
Kumar, following Slodowy \cite{Sl}, considers the basic representation
${\cal{H}}_{\ell}$ for a fixed
$\ell$, and a highest weight vector $v_{\ell}$. The subgroup $\LGp$ is the
stabilizer of the line
$kv_\ell$ in
$\proj({\cal{H}}_{\ell})$, so the map $g\mapsto gv_\ell$ induces  an injection
$i_\ell:\LG/\LGp\mono\proj({\cal{H}}_{\ell})$. Let $U$ be the subgroup of
$\LGp$ consisting of elements $A(z)$ such that $A(0)$ is in the unipotent part
of a fixed Borel subgroup $B\subset G$; to each element $w$ of the Weyl group
is
associated a ``Schubert variety"
$X_w$ which is a finite union of orbits of $U$. It turns out that the image
under
$i_\ell$ of  $X_w$
 is actually contained in some finite-dimensional projective subspace $\proj_w$
of
$\proj({\cal{H}}_{\ell})$, and is Zariski closed in $\proj_w$.   This defines
on
$X_w$ a structure  of reduced projective variety, and a structure of
ind-variety on
$\LG/\LGp=\limind X_w$.

By a result due to Faltings (\cf the Appendix of \cite{BL1} for $SL_{r}$), the
irreducible integrable representation
${\cal{H}}_{\ell}$  of $\Lgh$ can be ``integrated" to an {\it algebraic}
projective representation of $\LG$, that is a morphism of $k$-groups
$\LG\rightarrow PGL({\cal{H}}_{\ell})$.  It follows that the map $i_c$ is a
morphism of ind-schemes of $\Q$ into
$\proj({\cal{H}}_{\ell})$ (which is the direct limit of its finite-dimensional
subspaces).  But $i_\ell$ is even an {\it embedding}. It is injective by what
we said above; let us check that it induces an injective map on the tangent
spaces. Since it is equivariant under the action of
$\LG$ it is enough to prove this at  the origin $\omega$ of $\Q$. Then it
follows from the fact  that the annihilator of $v_\ell$ in the Lie algebra
$\Lg$ is
$\Lgp$.

Therefore the restriction of $i_\ell$ to each of the subvarieties  $\QN$ is
proper, injective, and injective on the tangent spaces, hence is an embedding
(in some finite-dimensional projective subspace of
$\proj({\cal{H}}_{\ell})$). Each $X_w$ is contained in some $\QN$, and
therefore is a closed subvariety of $\QNred$. Each orbit of $U$ is
contained in some $X_w$; since the
$X_w$'s define an ind-structure, each
$\QN$ is contained in some $X_w$, so that $\QNred$ is a
subvariety of
$X_w$. Since $\Q$ is the direct limit of the
$\QNred$, the two ind-structures coincide.
\end{proof}

\section{The ind-group $\LGX$}

Throughout this section  we suppose $k=\comp$ and $G$ semi-simple and simply
connected.

\begin{prop}\label{LGX-is-integral} The ind-group $\LGX$ is integral.
\end{prop}

\begin{cor}\label{LGX-has-no-characters} Every character $\chi:\LGX\ra G_m$ is
trivial.
\end{cor}

\begin{proof} The differential of $\chi$, considered as a function on $\LGX$,
is
everywhere vanishing. Indeed, since $\chi$ is a group morphism, this means that
the deduced Lie algebra morphism ${\goth{g}}\otimes A_{X}\ra k$ is zero. But as
the derived algebra  ${\cal D}({\goth{g}}\otimes A_X)$ is ${\cal
D}({\goth{g}})\otimes A_X$ and therefore equal to ${\goth{g}}\otimes A_X$
because
${\goth{g}}$ is simple, any Lie algebra morphism ${\goth{g}}\otimes A_X\ra k$
is
trivial.

As  $\LGX$ is integral we can write $\LGX$ as the direct limit of integral
varieties
$V_n$. The restriction of $\chi$ to $V_{n}$ has again zero derivative and is
therefore constant. For large $n$, the varieties $V_{n}$ contain $1$. This
implies
$\chi_{\mid V_{n}}=1$ and we are done.
\end{proof}

\begin{proof} To see that the ind-group $\LGX$ is reduced, consider the \'etale
trivial morphism $\bar\pi:{\cal{Q}}\ra\M$. Locally for the \'etale topology,
$\bar\pi$ is a product $\Omega\times\LGX$ (where
$\Omega$ is an \'etale neighborhood of
$\M$). Then use that ${\cal{Q}}$ is reduced by
(\ref{ind-structures-are-the-same}). As connected ind-groups are irreducible
by Proposition 3 of \cite{Sh} it is enough to show that $\LGX$ is connected.

The idea how to prove that $\LGX$ is connected is due to V. Drinfeld:
consider distinct points
$p_{1},\dots, p_{i}$ of $X$ which are all distinct from $p$. Define
$X^{*}_i=X\moins\{p,p_{1},\dots,p_{i}\}$ and, for every
$k$-algebra
$R$, define  $X_{i,R}^{*} = X^{*}_i\times_k\Spec(R)$. Denote by $A_{X_{i,R}}$
the
$k$-algebra
$\Gamma(X_{i,R}^{*},{\cal{O}}_{X^{*}_{i,R}})$ and by
 $\LGXi$ the $k$-group $R\mapsto G(A_{X_{i,R}})$. As $\LGX$, the $k$-group
$\LGXi$ is an ind-group (\cf \ref{Loop-groups}). The natural inclusion
$A_{X_{i,R}}\subset A_{X_{i+1,R}}$ defines a closed immersion
$f:\LGXi\ra\LGXip$.

\begin{lem}\label{pi_0} The closed immersion
$\LGXi\ra\LGXip$ defines a bijection
$$\pi_{0}(\LGXi)\isom\pi_{0}(\LGXip).$$
\end{lem}

\begin{proof} Consider the morphism $\LGXip\ra\LG$ defined by the developpement
in Laurent series at $p_{i+1}$.  We get a morphism
$\phi_{i+1}:\LGXip\ra\Qp$, where we denote $\Qp=\LG/\LGp$. (of course $\Q=\Qp$
but we emphasize here that we will consider the point
$p_{i+1}$ and not $p$.)

\medskip
\noindent{\em Claim:} The morphism $\phi_{i+1}:\LGXip\ra\Qp$ induces an
isomorphism on the level of stacks $\bar\phi_{i+1}:\LGXip/\LGXi\simeq\Qp$ and
is locally trivial for the
\'etale topology.
\medskip

The lemma reduces to the claim. Indeed, as $G$ is semi-simple and simply
connected, we have $\pi_{i}([\Qp]^{an})=1$ for $i=0,1$ (by
(\ref{ind-structures-are-the-same}) and Kumar and Mathieu) and the exact
homotopy sequence associated of the (Serre)-fibration $\phi_{i+1}$ shows that
$\pi_{0}([\LGXi]^{an})\isom\pi_{0}([\LGXip]^{an})$ ($[]^{an}$ means we
consider the usual topology). From the bijection
$\pi_{0}(\LGXiN)\isom\pi_{0}([\LGXiN]^{an})$ and Proposition 2 of \cite{Sh} it
follows then that
$\pi_{0}(\LGXi)\ra\pi_{0}([\LGXi]^{an})$ is bijective.

\medskip
\noindent{\em Proof of the claim:} Clearly $\phi_{i+1}:\LGXip\ra\Qp$ is $\LGXi$
invariant, hence defines a map $\bar\phi_{i+1}:\LGXip/\LGXi\ra\Qp$. Define a
morphism
$\Qp\ra\LGXip/\LGXi$ as follows. Let $R$ be a
$k$-algebra. By Proposition \ref{pairs} to an element of $\Qp(R)$ corresponds a
$G$-bundle
$E\ra X_R$ together with a trivialization $\tau_{p_{i+1}}:G\times
X_{p_{i+1},R}^*\ra E_{\mid X_{p_{i+1},R}^*}$. Here by $X_{p_{i+1},R}^*$ we
denote
$(X\moins\{p_{i+1}\})\times_k\Spec(R)$. For any $R$-algebra $R^\prime$, denote
$T(R^\prime)$ the set of trivializations $\tau_i$ of $E_{R^\prime}$ over
$X_{i,R}^{*}$. This defines a $R$-space $T$ on which $\LGXi$ acts. By Theorem
\ref{Drinfeld-Simpson} it is a torsor under $\LGXi$. For any $\tau_i\in
T(R^\prime)$
 the composite  $\tau_i^{-1}\circ\tau_{p_{i+1}}$ defines a morphism
 $X_{i+1,R}^{*}\ra G$ hence an element of $\LGXip(R)$. In this way we associate
 functorially to an object $(E,\tau_{p_{i+1}})$ of $\Qp(R)$  a
$\LGXi$-invariant map
$\alpha:T\ra\LGXip$, which defines the inverse of $\bar\phi$. The assertion
concerning the local triviality is proved as in Theorem
\ref{th:Uniformization}.
\end{proof}

Let us show that every element $g\in\LGX(k)$ is in the connected component of
the unit of $\LGX(k)$.  Using the well known fact that $G(K)$ is generated by
the standard unipotent subgroups
$U_{\alpha}(K)$, $\alpha\in\Delta$, we may suppose that $g$ is of the form
$\prod_{j\in J}\exp(f_{j}n_{j})$ where the $n_{j}$ are nilpotent elements of
$\g$ and
$f_{j}\in K$. Let
$\{p_{1},\dots,p_{i}\}$ be the poles of the functions $f_{j}, j\in J$. The
morphism
$$\begin{diagram} {\tr{A}}^{1}&\lra&\LGXi\\ t&\mapsto&\prod_{j\in
J}\exp(tf_{j}n_{j})\\
\end{diagram}
$$ is a path from $g$ to $1$ in $\LGXi$. By Lemma \ref{pi_0}, the morphism
$\pi_{0}(\LGX)\ra\pi_{0}(\LGXi)$ is bijective which proves that $g$ and $1$ are
indeed in the same connected component of $\LGX$.
\end{proof}

\section{Pfaffians}

Let $k$ be an algebraically closed field of characteristic $\not=2$ and $S$ a
$k$-scheme.

\subsection{The Picard categories}

Let $A$ be $\reln$ or $\reln/2\reln$. Denote by ${\goth{L}}_{A}$ the groupoid
of
$A$-graded invertible ${\cal{O}}_{S}$-modules. The objects of ${\goth{L}}_{A}$
are pairs
$[L]=(L,a)$ of invertible
${\cal{O}}_{S}$-modules $L$ and locally constant functions $a:S\ra A$,
morphisms
$[f]:[L]\ra[M]$ are defined if $a=b$ and are isomorphisms
$f:L\ra M$ of
${\cal{O}}_{S}$-modules. Denote $\idbb_{A}$ the object $({\cal{O}}_{S},0)$. The
category ${\goth{L}}_{A}$ has tensor products, defined by
$[L]\otimes[M]=(L\otimes M,a+b).$ Given
$[L]$ and $[M]$ we have Koszul's symmetry isomorphism $\sigma_{_{[L],[M]}}:
[L]\otimes[M]\ra[M]\otimes[L]$ defined on local sections $\ell$ and
$m$ by $\sigma_{_{L,M}}(\ell\otimes m)=(-1)^{ab}m\otimes\ell$.

Denote $\det_{A}$ the functor from the category of coherent locally free
${\cal{O}}_{S}$-modules with isomorphisms defined by
$\det_{A}=(\Lambda^{max},\rang(V))$ and $\det_{A}(f)=\Lambda^{max}(f)$.

In the following we drop the subscript $A$ for $A=\reln$ and replace it by $2$
for
$A=\reln/2\reln$.

\subsection{Pfaffians}\label{Pfaffians-generalities}

Let $V$ be a coherent locally free ${\cal{O}}_{S}$-module of rank $2n$. Let
$\Pf:\Lambda^{2}V^{*}\ra\Lambda^{2n}V^{*}$ be the unique map that commutes with
base changes and such that if $(e^{*}_{1},\dots,e^{*}_{2n})$ is a local frame
of
$V^{*}$ and
$\alpha=\Sigma_{i<j}a_{ij}e^{*}_{i}\wedge e^{*}_{j}$, then
$$\Pf(\alpha)=\pf(a)e_{1}^{*}\wedge\dots\wedge e^{*}_{2n}$$ where $\pf(a)$ is
the pfaffian [Bourbaki, Alg\`ebre 9.5.2] of the alternating matrix
$a_{ij}=-a_{ji}$ for
$i>j$ and $a_{ii}=0$ for $i=1, \dots, 2n$.

Suppose $\alpha:V\ra V^{*}$ is skewsymmetric. View $\alpha$ as a section of
$\Lambda^{2}V^{*}$ and define the {\em pfaffian} of $\alpha$ as the section
$\Pf(\alpha):{\cal{O}}_{S}\ra\Lambda^{2n}V^{*}$. By [Bourbaki, Alg\`ebre 9.5.2]
we know that
\begin{equation}\label{square}
\begin{diagram} {\cal{O}}_{S}\otimes{\cal{O}}_{S}&\lra&{\cal{O}}_{S}\\
\sfl{\Pf(\alpha)\otimes\Pf(\alpha)}{}&\swfl{}{\det(\alpha)}\\
\Lambda^{2n}V^{*}\otimes\Lambda^{2n}V^{*}\\
\end{diagram}
\end{equation} commutes and that, if $u$ is an endomorphism of $V^{*}$, then
\begin{equation}\label{sbc}\begin{diagram}
{\cal{O}}_{S}&\efl{\Pf(\alpha)}{}&\Lambda^{2n}V^{*}\\
\sfl{\Pf(u\alpha u^{*})}{}&\swfl{}{\det(u)}\\
\Lambda^{2n}V^{*}\\
\end{diagram}
\end{equation} commutes.

\subsection{The pfaffian functor}

We consider the following category ${\cal{A}}={\cal{A}}^{\bullet}(S)$: objects
are complexes of locally free coherent
${\cal{O}}_{S}$-modules concentrated in degrees $0$ and $1$ of the form
$$0\lra E\efl{\alpha}{} E^{*}\lra 0$$ with $\alpha$ skewsymmetric. Morphisms
between two such complexes $E^{\bullet}$ and
$F^{\bullet}$ are morphisms of complexes
$f^{\bullet}:E^{\bullet}\lra F^\bullet$ such that $f^{\bullet*}[-1]$ is a
homotopy inverse of
$f^{\bullet}$, \ie $f^{\bullet*}[-1]\circ f^{\bullet}$ and
$f^{\bullet}\circ f^{\bullet*}[-1]$ are  homotopic to the identity.

Let $\pi:{\goth{L}}\ra{\goth{L}}_{2}$ be the projection functor,
$\Delta:{\goth{L}}_{2}\ra{\goth{L}}_{2}$ be the functor defined by
$[L]\mapsto[L]\otimes[L]$ and $[f]\mapsto[f]\otimes[f]$ and
$\Det:{\cal{A}}\ra{\goth{L}}$ be the determinant functor \cite{KM} .

\begin{prop}\label{pfaffianfunctor} There is a natural functor,
$\Pf:{\cal{A}}^{\bullet}\ra{\goth{L}}_{2}$,  commuting with base changes, and a
natural isomorphism of functors:
$$\pi\circ\Det\isom\Delta\circ\Pf.$$ Moreover, if $f^{\bullet}:E^{\bullet}\lra
E^\bullet$ is homotopic to the identity then $\Pf(f^{\bullet})=\id$.
\end{prop}

\begin{proof} Define $\Pf$ on the level of objects by
$\Pf(E^{\bullet})=\det_{2}(E)$. On the level of morphisms we do the following.

Let $f^{\bullet}=(f_{0},f_{1}):E^{\bullet}\lra F^{\bullet}$ be a morphism of
${\cal{A}}^{\bullet}$:
$$\begin{diagram} E&\efl{\alpha_{E}}{}&E^{*}\cr
\sfl{f_{0}}{}&&\sfl{f_{1}}{}\cr F&\efl{\alpha_{F}}{}&F^{*}\cr
\end{diagram}
$$  and consider the complex $C^{\bullet}_{f}$ (which is up to sign the cone of
$f^{\bullet}$)
\bigskip\bigskip
$$C_{f}^{\bullet}= 0\lra E\efl{\begin{pmatrix}\alpha_{E}\cr
-f_{0}\end{pmatrix}}{} E^{*}\oplus  F\efl{(f_{1}\ \alpha_{F})}{}F^{*}\lra 0$$
As $f^{\bullet}$ is a quasi-isomorphism, $C^{\bullet}_{f}$ is acyclic. By the
usual additivity property of determinants, we get a canonical isomorphism
$$d(f):\Lambda^{max}E\otimes\Lambda^{max}F^{*}\ra\Lambda^{max}E^{*}
\otimes\Lambda^{max}F.$$ Recall that this isomorphism is defined by taking a
section
$\begin{pmatrix}u\\ v\end{pmatrix}$ of $(f_{1}\ \alpha_{F})$ and calculating
the
determinant, which is independent of this choice, of the morphism
$$M(f)=\begin{pmatrix}\alpha&u\cr -f_{0}&v\end{pmatrix}\in
\Hom(E\oplus F^{*},E^{*}\oplus F)$$

\begin{lem}\label{compagnion} There is a skew-symmetric morphism
$\gamma_{f}\in\Hom(F^{*},F)$ such that
$\begin{pmatrix}f_{0}^{*}\\ \gamma\end{pmatrix}$ is a section of
$(f_{1}\ \alpha_{F})$.
\end{lem}

\begin{proof} As $f\circ f^{*}[-1]$ is homotopic to $\Id$ there is a morphism
$h$ such that
$f_{0}f_{1}^{*}-1=h\alpha_{F}$ and $f_{1}f_{0}^{*}-1=\alpha_{F} h$. Now define
$\gamma_{f}=\frac{h^{*}-h}{2}$.
\end{proof}

The pfaffian of the skew-symmetric morphism
$$M(f,\gamma_{f})=\begin{pmatrix}\alpha&f_{0}^{*}\cr -f_{0}&
\gamma_{f}\end{pmatrix}\in\Hom(E\oplus F^{*},E^{*}\oplus F)$$ defines a section
$\pf(M(f,\gamma_{f})):{\cal{O}}_{S}\ra
\Lambda^{max}E^{*}\otimes\Lambda^{max}F$.

\begin{lem} The section $\pf(M(f,\gamma_{f}))$ is independent of the choice of
$\gamma_{f}$.
\end{lem}

\begin{proof} Suppose $\gamma_{f}^{\prime}$ is another morphism satisfying
(\ref{compagnion}). Then there is $g\in\Hom(F^{*},E)$ such that $\alpha_{E}g=0$
and
$f_{0}g=-g^{*}f_{0}^{*}$ [use that $\gamma_{f}$ and $\gamma_{f}^{\prime}$ are
skew]. These relations give
$$M(f,\gamma_{f}^{\prime})=
\begin{pmatrix}1&0\cr {\frac{g}{2}}^{*}&1\end{pmatrix}
M(f,\gamma_{f})\begin{pmatrix}1&\frac{g}{2}\cr 0&1\end{pmatrix}$$ which in turn
implies the required equality by (\ref{sbc}).
\end{proof}

As $\rank(E)=\rank(F)\bmod 2$, we get the isomorphism in ${\goth{L}}_{2}$:
$$\pf(M(f)):\idbb_{2}\isom\det_{2}(E)^{*}\otimes\det_{2}(F).$$ Define the
pfaffian of
$f^{\bullet}$ by
$$\Pf(f^{\bullet}):\det_{2}(E)\efl{1\otimes\pf(M(f,\gamma_{f}))}{}
\det_{2}(E)\otimes\det_{2}(E)^{*}\otimes\det_{2}(F)\efl{\ev_{_{\det_{2}(E)}}}{}
\det_{2}(F)$$

\begin{lem} $\Pf:{\cal{A}}\ra{\goth{L}}_{2}$ defines a functor.
\end{lem}

\begin{proof} As $\pf(M(\Id,0))=1$, we have $\Pf(\Id)=\Id$. Let
$f^{\bullet}:E^{\bullet}\ra F^{\bullet}$ and $g:F^{\bullet}\ra G^{\bullet}$ be
two morphisms of ${\cal{A}}$. Then the following diagram is commutative
$$
\begin{diagram}
\idbb_{2}&\efl{\pf(M(f,\gamma_{f}))\otimes\pf(M(g,\gamma_{g}))}{}&
\det_{2}(E)^{*}\otimes\det_{2}(F)\otimes\det_{2}(F)^{*}\otimes\det_{2}(G)\\
\sfl{\Id}{}&&\sfl{}{1\otimes\ev_{\det_{2}(F)}\otimes 1}\\
\idbb_{2}&\efl{\pf(M(g\circ f,\gamma_{g\circ f}))}{}&
\det_{2}(E)^{*}\otimes\det_{2}(G)\\
\end{diagram}
$$ Indeed, remark that
$\gamma_{gf}=g_{0}\gamma_{f}g_{0}^{*}+\gamma_{g}$ satisfies (\ref{compagnion})
for
$g\circ f$ and make use of (\ref{sbc}) first with
$$
\begin{pmatrix}
\alpha_{E}&f_{0}^{*}&0&f_{0}^{*}g_{0}^{*}\\
-f_{0}&\gamma_{f}&1&\gamma_{f}g_{0}^{*}\\ 0&-1&0&0\\
-g_{0}f_{0}&g_{0}\gamma_{f}&0&\gamma_{gf}\\
\end{pmatrix} =
\begin{pmatrix} 1&0&0&0\\ 0&1&0&0\\ -f_{1}&-\alpha_{F}&1&0\\ 0&g_{0}&0&1\\
\end{pmatrix}
\begin{pmatrix}
\alpha_{E}&f_{0}^{*}&0&0\\ -f_{0}&\gamma_{f}&0&0\\ 0&0&\alpha_{F}&g_{0}^{*}\\
0&0&-g_{0}&\gamma_{g}\\
\end{pmatrix}
\begin{pmatrix} 1&0&-f_{1}^{*}&0\\ 0&1&\alpha_{F}&g_{0}^{*}\\ 0&0&1&0\\
0&0&0&1\\
\end{pmatrix},
$$ and then with
$$
\begin{pmatrix}
\alpha_{E}&0&0&f_{0}^{*}g_{0}^{*}\\ 0&\gamma_{f}&1&0\\ 0&-1&0&0\\
-g_{0}f_{0}&0&0&\gamma_{gf}\\
\end{pmatrix} =
\begin{pmatrix} 1&0&f_{0}^{*}&0\\ 0&1&0&0\\ 0&0&1&0\\ 0&0&g_{0}\gamma_{f}&1\\
\end{pmatrix}
\begin{pmatrix}
\alpha_{E}&f_{0}^{*}&0&f_{0}^{*}g_{0}^{*}\\
-f_{0}&\gamma_{f}&1&\gamma_{f}g_{0}^{*}\\ 0&-1&0&0\\
-g_{0}f_{0}&g_{0}\gamma_{f}&0&\gamma_{gf}\\
\end{pmatrix}
\begin{pmatrix} 1&0&0&0\\ 0&1&0&0\\ f_{0}&0&1&-\gamma_{f}g_{0}^{*}\\ 0&0&0&1\\
\end{pmatrix}
$$ The commutativity of the above diagram shows $\Pf(g\circ
f)=\Pf(g)\circ\Pf(f)$.
\end{proof}

The statement on the natural transformation follows from the definitions and
(\ref{square}). It remains to prove that if $f^{\bullet}:E^{\bullet}\ra
E^{\bullet}$ is homotopic to the identity, then $\Pf(f)=\Id$. Indeed, let
$h:E^{*}\ra E$ be such that
$f_{0}-h\alpha_{E}=1$ and $f_{1}-\alpha h=1$. Then $\gamma_{f}:=-h+f_{0}h^{*}$
satisfies (\ref{compagnion}) and the statement follows from (\ref{sbc}) and
$$
\begin{pmatrix}a_{E}&1\\-1&0\end{pmatrix}=
\begin{pmatrix}1&0\\ h&1\end{pmatrix}
\begin{pmatrix}a_{E}&f_{0}^{*}\\ -f_{0}^{*}&\gamma_{f}\end{pmatrix}
\begin{pmatrix}1&h^{*}\\ 0&1\end{pmatrix}$$ This completes the proof of
Proposition
\ref{pfaffianfunctor}.
\end{proof}

\section{Line bundles on $\M$.}

Let $k$ be an algebraically closed field of characteristic $\not=2$ in Sections
(\ref{the-pf-bundle})-(\ref{dualizing-sheaf}). Denote
$\Pic(\M)$ the group of isomorphism classes of line bundles on $\M$. (See
(\cite{BL1}, 3.7) for a discussion of line bundles over $k$-spaces and stacks).
We will construct special elements of
$\Pic(\M)$.

\subsection{The determinant line bundle.}\label{the-det-bundle}

We start with the well known case of $G=GL_{r}$:  let ${\cal{F}}$ be a family
of
vector bundle of rank
$r$ parameterized by the locally noetherian $k$-scheme $S$.  Recall that the
complex
$Rpr_{1*}({\cal{F}})$ may be represented by  a perfect complex of length one
$K^{\bullet}$ and define ${\scr{D}}_{{\cal{F}}}$ to be
$\det(K^{\bullet})^{-1}$. This does not, up to canonical isomorphism, depend on
the choice of $K^{\bullet}$. As the formation of the determinant commutes with
base change, the fiber of
${\scr{D}}_{{\cal{F}}}$ over the point $s\in S$ is
$\Lambda^{max}H^{0}(X,{\cal{F}}(s))^{*}\otimes
\Lambda^{max}H^{1}(X,{\cal{F}}(s))$. The line bundle ${\scr{D}}_{{\cal{F}}}$ is
called the {\em determinant of cohomology} line bundle associated to the family
${\cal{F}}$.

Let ${\cal{U}}$ be the universal vector bundle on
$\MGL\times X$ and define the determinant line bundle
${\scr{D}}=\det(Rpr_{1*}{\cal{U}})^{-1}$. It has the following universal
property: for every family ${\cal{F}}$ of vector bundle parameterized by the
locally noetherian
$k$-scheme
$S$, we have
$f_{{\cal{F}}}^{*}({\scr{D}})={\scr{D}}_{{\cal{F}}}$ where
$f_{{\cal{F}}}:S\ra\MGL$ is the deduced modular morphism.

For the case of general $G$, consider a representation $\rho:G\ra GL_{r}$ and
consider the morphism obtained by extension of structure group
$f_{\rho}:\M\ra\MGL$. Then define the determinant of cohomology associated to
$\rho$ by ${\scr{D}}_{\rho}=f_{\rho}^{*}({\scr{D}})$.

\subsection{The Pfaffian bundle}\label{the-pf-bundle}

Consider $G=\Spin_{r}$ with $r\geq 3$ (resp. $G=G_2$). Then the standard
representation $\varpi_{1}$ factors through $SO_{r}$ (resp. $SO_{7}$). The
stack
$\MSO$ has two components: $\MSOzero$ and
$\MSOone$. They are distinguished by the second
Stiefel-Whitney class
$$w_{2}: H^{1}_{\acute
et}(X,SO_{r})\ra\H^{2}_{\acute et}(X,\reln/2\reln)=\reln/2\reln.$$

Let $\kappa$ be a theta-characteristic on $X$. Twisting by $\kappa$,  we may
and will see a $SO_{r}$-bundle as a vector bundle $F$ with trivial determinant
together with a {\em symmetric} isomorphism $\sigma:F\ra F^{\vee}$, where
$F^{\vee}=\ul\Hom_{{\cal{O}}_{X}}(F,\omega_{_{X}})$. The following Proposition
shows the existence, for every $\kappa$, of a canonical square root
${\scr{P}}_{\kappa}$ of the determinant bundle ${\scr{D}}_{\varpi_{1}}$ over
$\MSO$.

\begin{prop} Let $(F,\sigma)$ be a family of vector bundles $F$ equipped with a
quadratic form $\sigma$ with values in $\omega_{_{X}}$ parameterized by the
locally noetherian
$k$-scheme $S$. Then the determinant of cohomology ${\cal{D}}_{F}$ admits a
canonical square root ${\cal{P}}_{(F,\sigma)}$. Moreover, if $f:S^{\prime}\ra
S$ is a morphism of locally noetherian $k$-schemes then we have
${\cal{P}}_{(f^{*}F,f^{*}\sigma)}=f^{*}{\cal{P}}_{(F,\sigma)}$.
\end{prop}

\begin{proof} By (\cite{So2}, prop. 2.1 and proof of corollary 2.2, \cf also
\cite{Ke}), Zariski locally on $S$, there are length
$1$ complexes
$M^{\bullet}$ of finite free ${\cal{O}}_{S}$-modules and quasi-isomorphisms
$f:M^{\bullet}\ra Rpr_{1*}(F)$ such that the composition in the derived
category
$D(S)$ (use $\sigma$ and Grothendieck duality)
$$M^{\bullet}\efl{f}{} Rpr_{1*}(F)\efl{\tau}{}
\R\ul\Hom^{\bullet}( Rpr_{1*}(F),{\cal{O}}_{S})[-1]\efl{f^{*}[-1]}{}
M^{\bullet*}[-1]$$ lifts to a symmetric isomorphism of complexes
$\varphi:M^{\bullet}\ra M^{\bullet*}[-1]$:
$$
\begin{diagram} 0&\lra&M^{0}&\efl{d_{M^{\bullet}}}{}&M^{1}&\lra&0\\
&&\sfl{\varphi_{0}}{\wr}&\sefl{}{\alpha}&\sfl{\wr}{\varphi_{0}^{*}}\\
0&\lra&M^{1*}&\efl{-d_{M^{\bullet *}}}{}&M^{0*}&\lra&0\\
\end{diagram}
$$

\comment{ [Indeed, for every $s\in S$ we may find, in a Zariski neighborhood
$U$ of
$s$, a length $1$ complex $(M^{\bullet},d_{M^{\bullet}})$ of finite free
${\cal{O}}_{S}$-modules such that $d_{M^{\bullet}}(s)=0$ and a
quasi-isomorphism
(over $U$)
$f:M^{\bullet}\ra Rpr_{1*}(F)$. The isomorphism (in $D^{b}(U)$)
$f^{*}[-1]\circ\tau\circ f:M^{\bullet}\ra M^{\bullet*}[-1]$ lifts to a morphism
of complexes $\phi$, as the components of $M^{\bullet}$ are free. As $\tau$ is
symmetric, the symmetrization  $\varphi$ of $\phi$ still lifts
$f^{*}[-1]\circ\tau\circ f$ and as $\varphi$ is an isomorphism at the point
$s$,
$\varphi$ remains an isomorphism, after eventually shrinking $U$, over $U$.] }

Define $\widetilde{M}^{\bullet}$ by
$0\ra M^{0}\efl{\alpha}{}M^{0*}\ra 0.$ Then $\alpha$ is skew and we have a
natural isomorphism of complexes $\psi:M^{\bullet}\ra\widetilde{M}^{\bullet}$
such that
$\psi^{*}[-1]\psi=f^{*}[-1]\tau f$ in $D^{b}(S)$.

Cover $S$ by open subsets $U_{i}$ together with complexes
$(M_{i}^{\bullet},d_{M^{\bullet}_{i}})$ and quasi-isomorphisms
$f_{i}:M_{i}^{\bullet}\ra Rpr_{1*}(F)\restriction{U_{i}}$ as above.  We define
${\cal{P}}_{(F,\sigma)}$ over $U_{i}$ by
${\cal{P}}_{i,(F,\sigma)}=\Pf(\widetilde{M}^{\bullet}_{i})$ and construct
patching data
$\rho_{ij}:{\cal{P}}_{i,(F,\sigma)}\isom{\cal{P}}_{j,(F,\sigma)}$  over
$U_{ij}=U_{i}\cap U_{j}$ in the following way.
\comment{ Consider the diagram of isomorphisms in $D^{b}(U_{ij})$, with
$K_{ij}^{\bullet}=Rpr_{1*}(F)\restriction{U_{ij}}$
$$\begin{diagram}\widetilde{M_i}&\efl{{\psi_i^{-1}}}{}&M_{i}^{\bullet}&\\
&&\sfl{f_{i}}{}\cr M^{\bullet}_{j}&\efl{f_{j}}{}&K_{ij}^{\bullet}&\efl{\tau}{}&
K_{ij}^{\bullet*}[-1]&
\efl{f_{j}^{*}[-1]}{}& M^{\bullet*}_{j}[-1]\\
&&\sfl{\tau}{}&&&&\sfl{\wr}{\psi_j^*[-1]}\cr &&K_{ij}^{\bullet*}[-1]
&&&&\widetilde{M_j^\bullet}{}^*[-1]\cr &&\sfl{f_{i}^{*}[-1]}{}\cr
&&M^{\bullet*}[-1]_{i}\cr
\end{diagram}
$$ } Define first the morphism of complexes
$\Sigma_{ij}:\ \widetilde{M}_{i}\ra\widetilde{M}_{j}$ as a lifting of the
isomorphism in $D^{b}(U_{ij})$
$$\psi_{j}^{*-1}[-1]f_{j}^{*}[-1]\tau f_{i}\psi_{i}^{-1},$$  then $\rho_{ij}$
by
$\Pf(\Sigma_{ij})$ [note that it follows from the symmetry of $\sigma$ (and
that the components of the $\widetilde{M}_{i}$ are free) that $\Sigma_{ij}$ is
a morphism of
${\cal{A}}(U_{ij})$]. By \ref{pfaffianfunctor},
$\rho_{ij}$ does not depend on the particular chosen lifting and the
functoriality of $\Pf$ translates into
$\rho_{ii}=\Id$, $\rho_{ij}=\rho_{ik}\rho_{kj}$ and also
$\rho_{ij}=\rho_{ji}^{-1}$. Over $U_{i}$ we have
${\cal{P}}_{i,(F,\sigma)}\otimes{\cal{P}}_{i,(F,\sigma)}
=\det(\widetilde{M}_{i}^{\bullet})$. As usual, the
$\det(\widetilde{M}_{i}^{\bullet})$ path together (via $f_{i}\psi_{i}^{-1}$),
to
${\scr{D}}_{F}$ and we get, again by Proposition \ref{pfaffianfunctor}, a
canonical isomorphism
${\cal{P}}_{(F,\sigma)}\otimes {\cal{P}}_{(F,\sigma)}\isom {\cal{D}}_{F}.$
\end{proof}

\subsection{} Considering the universal family over $\MSO\times X$, we
get, by the above, for every theta-characteristic $\kappa$ a line bundle
${\scr{P}}_{\kappa}$ over $\MSO$. Consider
$$e:\MSpin\lra\MSOzero$$ defined by extension
of the structure group. This morphism defines a morphism on the level of
Picard groups hence we can define a line bundle, denoted by ${\scr{P}}$, which
is the pullback of the pfaffian line bundle ${\scr{P}}_{\kappa}$. We omit here
the index $\kappa$ as we will see that ${\scr{P}}$ on
$\MSpin$ does not depend on the choice of a particular
theta-characteristic. In the same way, we define the line bundle ${\scr{P}}$ on
$\MGtwo$.

\subsection{The pfaffian divisor.} Let $r\geq 3$ and $({\cal E},q)$ be the
universal quadratic bundle over
$\MSOzero\times X$. For $\kappa$ a theta-characteristic, let us
denote by
$\Theta_\kappa$ the substack defined by
$$\Theta_\kappa={\rm div}(Rpr_{1*}({\cal E}\otimes pr_2^*\kappa)).$$

{\em Claim: This substack is a divisor if and only if $r$ or $\kappa$ are
even.}

\begin{proof} Let $P=(E,q)$ be a $SO_r$-bundle, $r\geq 3$ and $\kappa$ be a
theta-characteristic. Then
\begin{equation}\label{w2} w_2(P)=\h^0(E\otimes\kappa)+r\h^0(\kappa)\bmod 2.
\end{equation} Indeed, by Riemanns invariance mod $2$ theorem, the right hand
side  of (\ref{w2}) denoted $\bar w_2(P)$ in the following, is constant over
the
$2$ connected components of $\MSO$. Because (\ref{w2}) is true at
the trivial $SO_{r}$-bundle ${\cal{T}}$, it is  enough to prove that $\bar w_2$
is not constant. Let $L,M\in J_{2}$ (where $J_{2}=$ points of order 2 of the
jacobian) such that for the Weyl pairing $<L,M>=1$. The choice of a
trivialization of their square defines a non degenerated quadratic form on
$E=(L\otimes M)\oplus L\oplus M\oplus (r-3){\cal{O}}_{X}$ hence a
$SO_r$-bundle $P$. By \cite{Mu}, we know that we have
$\bar w_2(P)=<L,M>\not=0=\bar w_2({\cal{T}})$, which proves (\ref{w2}). Choose
an ineffective theta-characteristic $\kappa_0$  and set
$L=\kappa_0\otimes\kappa^{-1}$.  If $r$ is even, there exists a $SO_r$-bundle
$P=(E,q)$ such that
$H^0(E\otimes\kappa)=0$ and $w_2(P)=0$ (choose $E=rL$ with $L\in J_{2}$ and use
(\ref{w2})). If $r$ is odd and $\kappa$ is even, there exists a $SO_r$-bundle
$P=(E,q)$ such that $H^0(E\otimes\kappa)=0$ and $w_2(P)=0$ (by (\cite{Be},
lemme
1.5), there is a $SL_2$-bundle $F$ on $X$ such that
$H^0(X,{\rm ad}(F)\otimes \kappa)=0$, then choose $E={\rm ad}(F)\oplus (r-3)L$
with the obvious  quadratic form.) If $r$ and $\kappa$ are odd, then
$H^0(E\otimes\kappa)$ is odd for all $P\in\MSOzero$.
\end{proof}

As the perfect complex $Rpr_{1*}({\cal E}\otimes pr_2^*\kappa)$ can
be locally represented by a skew-symmetric perfect complex of length one
$L\efl{\alpha}{}L^*$, the pfaffian of $\alpha$ defines  a local equation of an
effective divisor ${\Theta_\kappa\over 2}$ such that
$2{\Theta_\kappa\over 2}=\Theta_\kappa$. This gives an easier way to define, by
smoothness of $\M$, the pfaffian line bundle. The reason which motivated our
construction above was to define this square root for arbitrary quadratic
bundles (not only the even ones) and to make a construction for all
theta-characteristics and not only the even ones (when $r$ is odd).

\subsection{Proof of \ref{sr}.}\label{dualizing-sheaf} The dualizing line
bundle $\omega_{_{\M}}$ on $\M$ is by definition the determinant line bundle of
the cotangent complex of $\M$. Let $\Ad:G\ra\GL({\goth{g}})$ be the adjoint
representation. Then $\omega_{_{\M}}={\scr{D}}_{\Ad}^{-1}$.  Suppose that $G$
is semi-simple. Then the adjoint representation factors through the special
orthogonal group because of the existence of the Cartan-Killing form. Choose a
theta-characteristic
$\kappa$ on
$X$. Then, as in (\ref{the-pf-bundle}), we can define a square root
$\omega^{\frac{1}{2}}_{_{\M}}(\kappa)$ of
$\omega_{_{\M}}$.

\section{The Picard group of $\M$.}\label{section-pic-of-M}

Throughout the section suppose that $k=\comp$  and that $G$ is  simple and
simply connected.

Let $\Pic_{\LGX}(\Q)$ the group of $\LGX$-linearized line bundles on $\Q$.
Recall that a
$\LGX$-linearization of the line bundle ${\scr{L}}$ on $\Q$ is an isomorphism
$m^{*}{\scr{L}}\isom pr_{2}^{*}{\scr{L}}$, where
$m:\LGX\times\Q\ra\Q$ is the action of $\LGX$ on $\Q$, satisfying the usual
cocycle condition.

Consider the projection $\pi:\Q\ra\M$ of Theorem $\ref{th:Uniformization}$. Let
${\cal{L}}$ be a line bundle on $\M$. As $\pi^{*}$ induces an isomorphism
between the sections of ${\cal{L}}$ and
$\LGX$-invariant sections of $\pi^{*}{\cal{L}}$ (\cite{BL1}, Lemma 7.2), we
have

\begin{prop}\label{pullback-is-injective} The projection $\pi:\Q\ra\M$ induces
an injection
$$\pi^{*}:\Pic(\M)\hookrightarrow\Pic_{\LGX}(\Q).$$
\end{prop}

Any $\LGX$-linearization is necessarily unique:

\begin{prop}\label{unique-linearization} The forgetful morphism
$\Pic_{\LGX}(\Q)\ra\Pic(\Q)$ is injective.
\end{prop}

\begin{proof} The kernel of this morphism consists of the
$\LGX$-linearizations of the trivial bundle. Any to such trivializations differ
by an automorphisms of $pr_{2}^{*}{\cal{O}}_{\Q}$ that is by an invertible
function  on
$\LGX\times\Q$. Since $\Q$ is integral (\ref{ind-structures-are-the-same}), it
is the direct limit of the  integral projective varieties
$\Q=\limind\QNred$ and this function is the pull back of an
invertible function $f$ on
$\LGX$. The cocycle conditions on the linearizations imply that
$f$ is a character, hence $f=1$ by Lemma \ref{LGX-has-no-characters}.
\end{proof}

\subsection{$\Pic(Q(G))$ and the canonical central extension of $LG$.}
\label{Pic(Q)}

Consider the embedding
$\LG/\LGp\hookrightarrow\proj({\cal{H}}_1)$ of
(\ref{ind-structures-are-the-same}) and define ${\cal{O}}_{\Q}(1)$ as the
pullback of
${\cal{O}}_{\proj({\cal{H}}_1)}(1)$. By \cite{Ma} and \cite{KNR}, we know
$\Pic(\Q)=\reln{\cal{O}}_{\Q}(1)$. The $k$-group $\LG$ acts on $\Q$ but the
action does not lift to an action of $\LG$ on ${\cal{O}}_{\Q}(1)$. There is a
canonical device to produce an extension $\LGh$ of $\LG$ such that the induced
action of $\LGh$ lifts to an action of $\LGh$ on ${\cal{O}}_{\Q}(1)$: the
Mumford group. This is the group of pairs $(g,f)$ where $g\in\LG$ and
$f:g^{*}{\cal{O}}_{\Q}(1)\isom{\cal{O}}_{\Q}(1)$. We get a central extension
(note that $\Q$ is direct limit of projective integral schemes)
\begin{equation}\label{can-ext} 1\lra G_{m}\lra\LGh\lra\LG\lra 1
\end{equation} Note that the Mumford group of
${\cal{O}}_{\proj({\cal{H}}_1)}(1)$ is
$GL({\cal{H}}_1)$. As the projective representation ${\cal{H}}_1$ of $\Lg$ can
be integrated to a projective representation $\phi:\LG\ra PGL({\cal{H}}_1)$
(\cf
\ref{ind-structures-are-the-same}), by functoriality of the Mumford group,
(\ref{can-ext}) is also the pullback by $\phi$ of the central extension
\begin{equation}\label{gl(H)-ext} 1\lra\G_{m}\lra GL({\cal{H}}_1)\lra
PGL({\cal{H}}_1)\lra 1.
\end{equation} Moreover, the restriction of
${\cal{H}}_1$ to $\Lgp$ can be integrated to a representation
$\LGp\ra GL({\cal{H}}_{1})$. It follows that (\ref{gl(H)-ext}) splits {\em
canonically} over $\LGp$, and ${\cal{O}}_{\Q}(1)$ is  the line bundle on the
homogeneous space $\Q=\widehat\LG/\widehat\LGp$ associated to the character
$G_m\times\LGp\ra G_m$ defined by the first projection.

\subsection{}\label{central-extensions-and-Dynkin-index} Consider the case of
$G=SL_{r}$. In this case, we know from \cite{BL1}, that the pullback of the
determinant line ${\cal{D}}$ bundle is ${\cal{O}}_{\QSL}(1)$.  It
follows from (\ref{pullback-is-injective}), (\ref{unique-linearization}) and
(\ref{Pic(Q)}) that $\Pic(\M)=\reln{\scr{D}}$. If $\rho:G\ra SL_{r}$ is a
representation of $G$ we get a commutative diagram
$$\begin{diagram}
\Q&\efl{}{}&\QSL\\
\sfl{}{}&\sefl{}{\varphi}&\sfl{}{}\\
\M&\efl{}{}&\MSL\\
\end{diagram}
$$

\begin{lem}\label{Pullback-is-L-to-Dynkin} Let $d_\rho$ be the Dynkin index of
$\rho$. Then the pullback of the determinant bundle under $\varphi$ is
${\cal{O}}_{\Q}(d_{\rho})$
\end{lem}

\begin{proof} Consider the pullback diagram of (\ref{can-ext}) for
$\LSL$:
$$
\begin{diagram} 1&\lra& G_m&\lra&\widetilde\LG&\lra&\LG&\lra&1\\
&&\parallel&&\sfl{}{}&&\sfl{}{}\\ 1&\lra& G_m&\lra&\widehat{\LSL}&\lra&
\LSL&\lra& 1
\end{diagram}
$$ Looking at the differentials (note that $\Lie(\LSLh)=\Lslh$ by \cite{BL1}),
on the level of Lie algebra, we restrict the universal central extension of
$\Lsl$ to $\Lg$. The resulting extension $\widetilde\Lg$ is (\cf Section
(\ref{section-Dynkin-index})) the Lie algebra of the Mumford group of
${\cal{O}}_{\Q}(d_{\rho})$ where $d_{\rho}$ is the Dynkin index of $d_{\rho}$,
which proves the lemma.
\end{proof}

\begin{cor} As a pullback, the line bundle
${\cal{O}}_{\Q}(d_\rho)$ is  $\LGX$-linearized.
\end{cor}

\subsection{Proof of Theorem \ref{th:Pic}.} By the above for series $A$ and
$C$, as the Dynkin index for the standard representation is $1$, that all line
bundles on
$\Q$ are $\LGX$-linearized. For the series $B$ and $D$ (and also for $G_{2}$)
the Dynkin index of the standard representation is $2$. But by the existence of
the pfaffian line bundles we see also in this case, that all line bundles on
$\Q$ are
$\LGX$-linearized and Theorem \ref{th:Pic} for $n=0$ follows from
(\ref{pullback-is-injective}), (\ref{unique-linearization}) and (\ref{Pic(Q)})

\begin{rem}\label{unique-splitting} The restriction to $\LGX$ of the canonical
central extension $\LGh$ of $\LG$ splits (at least for classical $G$ and
$G_2$).
Moreover this splitting is unique.
\end{rem}

\begin{proof} As ${\cal{O}}_{\Q}(1)$ admits a $\LGX$-linearization, the action
of
$\LGX$ on $\Q$ induced by the embedding $\LGX\subset\LG$ lifts to an action to
${\cal{O}}_{\Q}(1)$. The central extension of $\LGX$ obtained by pullback of
the
canonical central extension $\LGh$ of $\LG$ is the Mumford group of $\LGX$
associated to ${\cal{O}}_{\Q}(1)$. But this extension splits as the action
lifts and we are done. Two splittings differ by a character of $\LGX$. As there
is only the trivial character (corollary \ref{LGX-has-no-characters}) the
splitting must be unique.
\end{proof}

\section{Parabolic $G$-bundles.}\label{Gpar}

Throughout this section $G$ is simple, simply connected and $k=\comp$. We will
extend the previous sections to the moduli stacks of parabolic
$G$-bundles.

\subsection{}\label{Preparation-on-G/P} We use the notations of Section
\ref{Lie-theory}. We will recall some standard facts for Lie groups, which we
will use later. Let $G$ be the simple and simply connected algebraic group
associated to
$\g$. Denote by $T\subset G$ the Cartan subgroup associated to
${\goth{h}}\subset\g$ and by $B\subset G$ the Borel subgroup associated to
${\goth{b}}\subset\g$.  Given a subset $\Sigma$ of the set of simple roots
$\Pi$
(nodes of the Dynkin diagram), we can define a subalgebra
${\goth{p}}_{\Sigma}={\goth{b}}\oplus(\osum_{\alpha\in\Sigma}
{\goth{g}}_{-\alpha})\subset\g,$ hence a subgroup $P_{\Sigma}\subset G$. Remark
that
$P_{\emptyset}=G/B$,
$P_{\Pi}=G$ and that all $P_\Sigma$ contain $B$. The subgroup $P_\Sigma$ is
parabolic and conversely  any {\em standard } (\ie containing $B$) parabolic
subgroup arises in this way. Fix $\Sigma\subset \Pi$ and let
$\Gamma=\Pi\moins\Sigma$. Denote by $X(P_{\Sigma})$ the character group of
$P_\Sigma$. Any weight $\lambda$ such that
$\lambda(H_{\alpha})=0$ for all $\alpha\in\Sigma$ defines, via the exponential
map, a character of $P_{\Sigma}$ and all characters arise in this way,
\ie
$X(P_{\Sigma})=\{\lambda\in P/\lambda(H_{\alpha})=0
\text{ for all }\alpha\in\Sigma\}.$  Given $\chi\in X(P_\Sigma)$ we can define
the line bundle
$L_\chi=G\times^{P_\Sigma}k_\chi$ on the homogeneous space $G/P_\Sigma$. In
general, there is an exact sequence (\cite{FI}, prop. 3.1)
$$1\lra X(G)\lra X(P_\Sigma)\lra\Pic(G/P_\Sigma)\lra\Pic(G)\lra\Pic(P)\lra 0.$$
As
$G$ is simple, we have $X(G)=0$ and as $G$ is simply connected, we have
$\Pic(G)=0$ (\cite{FI}, cor. 4.5). We get the isomorphism
$X(P_{\Sigma})\isom\Pic(G/P_{\Sigma}).$ In particular, the Picard group of
$G/P_{\Sigma}$ is isomorphic to the free abelian group generated over $\Gamma$.

\subsection{} Consider closed points $p_{1},\dots,p_{n}$ of $X$, labeled with
the standard parabolic subgroup $P_{1},\dots,P_{n}$.  Let
$\Sigma_{1},\dots,\Sigma_{n}$ be the associated subsets of simple roots and
$\Gamma_{i}=\Pi\moins \Sigma_{i}$ for $i\in\{1,\dots,n\}$. In the following,
underlining a character will mean that we consider the associated sequence,
e.g.
$\ul{P}$ will denote the sequence
$(P_{1},\dots,P_{n})$, {\em etc.} Let $E$ be a $G$-bundle. As $G$ acts on
$G/P_{i}$ we can define the associated
$G/P_{i}$-bundle $E(G/P_{i})$.

\begin{defi} (\cf \cite{MS}) A quasi-parabolic $G$-bundle of type $\ul{P}$ is a
$G$-bundle $E$ on $X$ together with, for all $i\in\{1,\dots,n\}$, an element
$F_{i}\in E(G/P_{i})(p_{i})$. A parabolic $G$-bundle of type $(\ul{P},\ul{m})$
is a quasi-parabolic $G$-bundle of type
$\ul{P}$ together with, for $i\in\{1,\dots,n\}$, parabolic weights
$(m_{i,j})_{j\in\Gamma_{i}}$ where the $m_{i,j}$ are strictly positive
integers.
\end{defi}

\subsection{} Let $R$ be a $k$-algebra, $S=\Spec(R)$. A family of
quasi-parabolic
$G$-bundles of type $\ul{P}$ parameterized by $S$ is a $G$-bundle $E$ over
$S\times X$ together with $n$ sections $\sigma_{i}:S\ra E(G/P_{i})_{\mid
S\times\{p_{i}\}}$. A morphism from $(E,\ul\sigma)$ to
$(E^{\prime},\ul\sigma^{\prime})$ is a morphism
$f:E\ra E^{\prime}$ of $G$-bundles such that for all $i\in\{1,\dots,n\}$ we
have
$\sigma^{\prime}_{i}=f_{\mid S\times\{p_{i}\}}\circ\sigma_{i}$.

We get a functor from the category of $k$-algebras to the category of groupoids
by associating to $R$ the groupoid having as objects families of
quasi-parabolic
$G$-bundles of type $\ul{P}$ parameterized by $S=\Spec(R)$ and as arrows
isomorphisms between such families. Moreover for any morphism $R\ra R'$ we have
a natural functor between the associated groupoids. This defines the $k$-stack
of quasi-parabolic $G$-bundles  of type
$\ul{P}$ which we will denote by $\Mpar$. The stack $\Mpar$ has, as $\M$, a
natural interpretation as a double quotient stack. Define
$$\Qpar=\Q\times\prod_{i=1}^{n} G/P_{i}.$$ The ind-group $\LGX$ acts on $\Q$
and by evaluation
$ev(p_{i}):\LGX\ra G$ at $p_{i}$ also on each factor $G/P_{i}$. We get a
natural
action of $\LGX$ on $\Qpar$. The analogue of Theorem
\ref{th:Uniformization}. for quasi-parabolic $G$-bundles is

\begin{th}\label{th:ParUniformization} (Uniformization) There is a canonical
isomorphism of stacks
$$\overline{\pi}:\LGX\bk\Qpar\isom\Mpar.$$ Moreover the projection map is
locally trivial for the \'etale topology.
\end{th}

\begin{proof} Let $R$ be a $k$-algebra, $S=\Spec(R)$. To an element
$(E,\rho,\ul{f})$ of $\Qpar(R)$ (with $f_{i}\in\Mor(S,G/P_{i})$), we can
associate a family of quasi-parabolic $G$-bundles of type $\ul{P}$
parameterized by $S$ in the following way. We only have to define the sections
$\sigma_{i}$:
$$\sigma_{i}:S\hfl{(\id,f_{i})}{} S\times G/P_{i}\hfl{\rho_i(G/P_{i})}{}
E(G/P_{i})_{\mid S\times\{p_{i}\}}.$$ We get a $\LGX$ equivariant map
$\pi:\Qpar\ra\Mpar$ which induces the map on the level of stacks
$\overline{\pi}:\LGX\bk\Qpar\ra\Mpar$.

Conversely, let $(E,\ul\sigma)$ be a family of quasi-parabolic $G$-bundles of
type
$\ul{P}$ parameterized by $S=\Spec(R)$. For any $R$-algebra $R^{\prime}$, let
$T(R^{\prime})$ be the set of trivializations $\rho$ of $E_{R^{\prime}}$ over
$X_{R^{\prime}}$. This defines a $R$-space $T$ which by Theorem
\ref{Drinfeld-Simpson} is a $\LGX$-torsor. To any element in $T(R^{\prime})$,
we can associate the family
$\ul{f}$ by
$$f_{i}:S\efl{\sigma_{i}}{} E(G/P_{i})_{\mid S\times\{p_{i}\}}
\efl{\rho_i(G/P_{i})^{-1}}{} S\times G/P_{i}\efl{pr_{2}}{} G/P_{i}.$$ In this
way we associate functorially to objects
$(E,\ul\sigma)$ of $\Mpar(R)$
$\LGX$-equivariant maps $\alpha:T\rightarrow \Qpar$. This defines a morphism of
stacks $$\Mpar\lra\LGX\bk\Qpar$$ which is the inverse of $\overline{\pi}$. The
second statement is clear from the proof of Theorem
\ref{th:Uniformization}.
\end{proof}

\subsection{} We study first line bundles over $\Qpar$. Using (\ref{Pic(Q)}),
(\ref{Preparation-on-G/P}) and $H^{1}(G/P_{i},{\cal{O}})=0$, we obtain the
following proposition, proving, as $\LGX$ has no characters, Theorem
\ref{th:Pic}.

\begin{prop} We have
$$\Pic(\Qpar)=\reln{\cal{O}}_{\Q}(1)\times\prod_{i=1}^{n}\Pic(G/P_{i})
=\reln{\cal{O}}_{\Q}(1)\times\prod_{i=1}^{n}X(P_{i}).$$
\end{prop}

\comment{Let $(E,\ul\sigma)$ be a family of quasi-parabolic
$G$-bundles of type $\ul{P}$ parameterized by the $k$-scheme $S=\Spec(R)$. Fix
$i\in\{1,\dots,n\}$ and $j\in\Gamma_{i}$. We may view $E\ra E(G/P_{i})$ as a
$P_{i}$-bundle. Therefore the character of $P_{i}$ defined by $-\varpi_{j}$
defines a line bundle on
$E(G/P_{i})$, hence by pullback, using the section
$\sigma_{i}:S\ra E(G/P_{i})_{\mid S\times\{p_{i}\}}$, a line bundle
${\scr{L}}_{i,j}$ over $S$. This works for any $S$ and we get a line bundle
over the stack $\Mpar$ which we denote again by ${\scr{L}}_{i,j}$.}

\section{Conformal blocs and generalized theta
functions.}\label{Identification}

Throughout this section $G$ is simple and simply connected and $k=\comp$.

\subsection{} Fix an integer $\ell\geq 0$ (the level) and let
$p_{1},\dots,p_{n}$ be distinct closed points of $X$ (we allow $n=0$ \ie no
points), each of it labeled with a dominant weight $\lambda_{i}$ lying in the
fundamental alc\^ove $P_{\ell}$.  Choose also another point $p\in X$, distinct
from the points $p_{1},\dots,p_{n}$. Define
$${\cal{H}}_{\ul\lambda}=
{\cal{H}}_{\ell}\otimes(\omal_{i=1}^{n}L_{\lambda_{i}}).$$ and let $\LgX$ be
$\g\otimes A_{X}$. We can map $\LgX$ via the Laurent developpement at the point
$p$ to $\Lg$. The restriction to
$\LgX$ of the universal central extension $\Lgh$ of $\Lg$  splits by the
residue
theorem, hence
$\LgX$ may be considered as a {\em sub Lie-algebra} of $\Lgh$. In particular,
${\cal{H}}_{\ell}$ is a $\LgX$-module. Evaluating $X\otimes f\in\LgX$ at the
point
$p_{i}$, we may consider $L_{\lambda_{i}}$ as a $\LgX$-module. Therefore
${\cal{H}}_{\ul\lambda}$ is a (left) $\LgX$-module. Define the space of
conformal blocks (or vacua) by
$$V_{X}(\ul{p},\ul\lambda)=[{\cal{H}}_{\ul\lambda}^{*}]^{\LgX}:=
\{\psi\in {\cal{H}}_{\ul\lambda}^{*}\ /\ \psi.(X\otimes f)=0\ \forall X\otimes
f\in\LgX\}.
$$ This definition is Beauville's description \cite{B} (see also \cite{So3}) of
the space of conformal blocks of  Tsuchiya, Ueno and Yamada \cite{TUY}.

The labeling of the points $p_{i}$ induces
$\Sigma_{i}=\{\alpha\in \Pi/\lambda_{i}(H_{\alpha})=0\}$,
$\Gamma_{i}=\Pi\moins\Sigma_{i}$ and $m_{i,j}=\lambda_{i}(H_{\alpha_{j}})$ for
$j\in\Gamma_{i}$, that is the type of a parabolic $G$-bundle. In particular we
get, for $\ell\in\natn$, a natural line bundle on the moduli stack $\Mpar$
defined by
$${\scr{L}}(\ell,\ul{m})={\scr{L}}^{\ell}\extern\bigl(\extern_{i=1}^{n}
(\extern_{j\in\Gamma_{i}}\reln{\scr{L}}_{i,j}^{m_{i,j}})\bigr).
$$ By construction, for the pull back of ${\scr{L}}(\ell,\ul{m})$ to $\Qpar$ we
have
$$\pi^{*}{\scr{L}}(\ell,\ul{m})={\cal{O}}_{\Q}(\ell)\extern\bigl(
\extern_{i=1}^{n}{\scr{L}}_{-\lambda_{i}}\bigr)$$ where
${\scr{L}}_{-\lambda_{i}}$ is the line bundle on the homogeneous space
$G/P_{i}$ defined by the character corresponding to the weight
$-\lambda_{i}$.

\subsection{Proof of (\ref{Verlinde}):}
We extend the method of \cite{BL1} and \cite{P}.

\medskip\noindent {\em Step 1:} As a pullback, $\pi^{*}{\scr{L}}(\ell,\ul{m})$
is canonically
$\LGX$-linearized, that is equipped with
$\varphi:m^{*}(\pi^{*}{\scr{L}}(\ell,\ul{m}))\isom
pr_{2}^{*}(\pi^{*}{\scr{L}}(\ell,\ul{m}))$. Denote by
$[H^{0}(\Qpar,\pi^{*}{\scr{L}}(\ell,\ul{m}))]^{\LGX}$ the space of
$\LGX$-invariant sections, that is the sections $s$ such that
$\varphi(m^{*}s)=pr_{2}^{*}s$. By Lemma 7.2 of \cite{BL1} we have the canonical
isomorphism
$$H^{0}(\Mpar,{\scr{L}}(\ell,\ul{m}))\isom
[H^{0}(\Qpar,\pi^{*}{\scr{L}}(\ell,\ul{m}))]^{\LGX}$$ Denote by
$[H^{0}(\Qpar,\pi^{*}{\scr{L}}(\ell,\ul{m}))]^{\LgX}$ the sections annihilated
by
$\Lie(\LGX)=\LgX$. By Proposition 7.4 of \cite{BL1}, using that $\LGX$ and
$\Qpar$ are integral (\ref{LGX-is-integral} and
\ref{ind-structures-are-the-same}), we have the canonical isomorphism
$$ [H^{0}(\Qpar,\pi^{*}{\scr{L}}(\ell,\ul{m}))]^{\LGX}\isom
[H^{0}(\Qpar,\pi^{*}{\scr{L}}(\ell,\ul{m}))]^{\LgX}
$$

\medskip\noindent {\em Step 2:} By definition of $\LGh$, the space
$H^{0}(\Qpar,\pi^{*}{\scr{L}}(\ell,\ul{m}))$ is naturally a $\LGh$-module.
Moreover we know that $\LGh$ splits over $\LGX$ (at least for classical $G$ and
$G_2$) and that this splitting is {\em unique}. The action of $\LgX\subset\Lgh$
deduced from this inclusion on
$H^{0}(\Qpar,\pi^{*}{\scr{L}}(\ell,\ul{m}))$ is therefore the same as the
preceding one.

\medskip\noindent {\em Step 3:} We have the canonical isomorphism of
$\Lgh$-modules
$$H^0(\Qpar,\pi^{*}({\scr{L}}(\ell,\ul{m})))\isom
H^{0}(\Q,{\cal{O}}_{\Q}(\ell))\otimes\bigl(\omal_{i=1}^{n}
H^{0}(G/P_{i},{\scr{L}}_{-\lambda_i})\bigr)$$ To see this apply the Kunneth
formula to the restriction of
${\scr{L}}(\ell,\ul{m})$ to the projective varieties
$\Qpar^{(N)}=\QN\times\prod_{i=1}^{n}G/P_{i}$, then use that inverse limits
commute with the tensor products by finite dimensional vector spaces.

\medskip\noindent {\em Step 4:} We have the canonical isomorphism of
$\LGh$-modules
$$H^{0}(\Q,{\cal{O}}_{Q}(\ell))\otimes\bigl(\omal_{i=1}^{n}
H^{0}(G/P_{i},{\scr{L}}_{-\lambda_i})\bigr)\isom {\cal{H}}_{\ell,0}^{*}\otimes
\bigl(\omal_{i=1}^{n}L_{\lambda_i}^{*}\bigr)$$ This is Borel-Bott-Weil theory,
in the version of Kumar-Mathieu (\cite{Ku},
\cite{Ma}) for the first factor, and the standard version \footnote{In
\cite{Bott} only the case G/B (\ie $\Sigma=\emptyset$) is considered but the
generalization to arbitrary
$G/P_{\Sigma}$ is immediate (and well known)} for the others.

The theorem follows from steps 1 to 4.

As we know the dimensions (at least for classical $G$ and $G_{2}$) for the
conformal blocks (\cite{F},\cite{B}, or \cite{So3} for an overview) we get the
Verlinde dimension formula for the spaces of generalized parabolic
theta-functions.

\comment{
\begin{cor}\label{cor:Verlinde-formula}(Verlinde formula) The dimension of the
space of generalized parabolic theta-functions
$H^0(\Mpar,{\scr{L}}(\ell,\ul{m}))$ is
$$(\# T_{\ell})^{g-1}
\sum_{\mu\in P_{\ell}}\Tr_{L_{\ul\lambda}}(\exp \frac{2\pi
i}{\ell+g^{*}}(\mu+\rho))
\prod_{\alpha\in\Delta_{+}}
\left| 2\sin \frac{\pi}{\ell+g^{*}}(\alpha,\mu+\rho)\right|^{2-2g}$$ with
$\#T_{\ell}=(\ell+g^{*})^{\rank{\goth{g}}}\#(P/Q)\#(Q/Q_{lg}).$
\end{cor} }

\section{Moduli spaces.}

\subsection{}\label{gen-on-CM} Suppose $char(k)=0$. We will show how the
previous results apply to the {\em coarse moduli spaces} of principal
$G$-bundles. We suppose that
 $G$ is reductive  and that $g\geq 2$. Recall  that a $G$-bundle $E$ over
$X$ is {\em semi-stable} (resp. {\em stable}) if for every parabolic subgroup
$P$ and for every reduction $E_{P}$ of $E$ to $G$, we have for every dominant
character (with respect to some Borel $B\subset P$)
$\chi$ of $P$, trivial over $Z_0(G)$, the following inequality
$\deg(E_P(\chi))\leq 0 \text{ (resp. $<$)}.$ A stable $G$-bundle $E$ is called
{\em regularly stable}, if moreover $\Aut(E)/Z(G)=\{1\}$.

Topologicially, $G$-bundles over $X$ are classified by  elements of
$\pi_{1}(G)$. By Ramanathan's \cite{Ra} theorem, there are coarse moduli spaces
$\Mt$ of semi-stable principal $G$-bundles of dimension
$(g-1)\dim G+\dim Z_{0}(G)$, which are irreducible, once the topological type
$\tau\in\pi_{1}(G)$ is fixed. Moreover $\Modt$ is normal and  the open subset
$\Modtreg\subset\Modt$ corresponding to  regularly stable $G$-bundles is
smooth.

\subsection{}\label{locally-factorial} Denote $\Cl$ the group of Weil divisor
classes. There is a commutative diagram
$$\begin{diagram}
\Pic(\Modt)&\efl{c}{}&\Cl(\Modt)\\
\sfl{r_{1}}{}&&\sfl{}{r_{2}}\\
\Pic(\Modtreg)&\efl{c_{reg}}{}&\Cl(\Modtreg)\\
\end{diagram}
$$ By normality, the restriction $r_{1}$ is injective, by smoothness of
$\Modtreg$, the canonical morphism $c_{reg}$ is an isomorphism and  as
(\cite{F1}, II.6)
$\codim_{\Modt}\Modt\setminus\Modtreg\geq 2$ (except when $g=2$ and $G$ maps
nontrivially to $PGL_2$) the restriction $r_{2}$ is an isomorphism. In
particular, $\Modt$ is locally factorial \cite{DN} if and only if $r_{1}$ is
surjective.

\subsection{}\label{det-for-CM} Consider $G=GL_{r}$. Then we may present
$\Modzero$ as the good quotient $\HG/GL(M)$ where $\HG$ is Grothendiecks
Quot scheme
$\HG=\Quot^{ss}(k^{M}\otimes{\cal{O}}_{X}(-N),P)$ parameterizing equivalence
classes (with the obvious equivalence relation) of pairs $[E,\alpha]$ with $E$
a
semi-stable vectorbundle of degree
$0$ and $\alpha:k^{M}\isom E(N)$, where $N$  and $M=rN+\chi(E)$. Let
${\cal{E}}$ be the universal family over
$\HG\times X$ and consider
$D=\det(Rpr_{1*}({\cal{E}}\otimes pr_{2}^{*}(L)),$ with $L$ a line bundle. It
is well known that $[E,\alpha]\in\HG$ has closed orbit exactly when
$E$ is polystable, \ie direct sum of stable bundles:
$E\simeq E_{1}^{\oplus n_{1}}\oplus\dots\oplus E_{\ell}^{\oplus n_{\ell}}$, and
that the action of the stabilizer $GL(n_{1})\times\dots GL(n_{\ell})$ is given
by the character
$$(g_1,\dots,g_\ell)
\mapsto\det(g_1)^{\chi(E_{1}\otimes L)}\cdot\dots\cdot
\det(g_\ell)^{\chi(E_{1}\otimes L)}.$$ Choose a line bundle $L$ of degree $g-1$
on $X$. Then $\chi(E_{q}\otimes L)=0$ for $q\in\HG$ and the action is
trivial. By Kempf's lemma \cite{DN}, $D$ descends to the determinant of
cohomology line bundle on $\ModGLzero$.

\subsection{Proof of \ref{th:Pic(Mmod)}.} Suppose $G$ is simple and simply
connected. We have (except for $g=2$ and $G=SL_2$)
$$\codim_{\M}(\M\setminus\Mreg)\geq 2.$$ To see this define
the Harder-Narasimhan filtration in the case of $G$-bundles and calculate the
codimension of the strata (\cite{LR}, Section 3) to show that for the open
substack $\Mss\subset\M$ corresponding to semi-stable $G$-bundles we  have
$\codim_{\M}(\M\setminus\Mss)\geq 2$, then use (\cite{F1},
II.6). The smoothness of $\M$ implies
$\Pic(\Mreg)=\Pic(\M)$ and it follows from Theorem \ref{th:Pic} that
$\Pic(\Mod)$ is an infinite cyclic group (note that the canonical morphism
$\Mreg\ra\Modreg$ induces an injection on the level of Picard groups). By
(\ref{det-for-CM}) and (\ref{section-Dynkin-index}), we know that the generator
is the determinant of cohomology for $G$ of type $A$ and $C$. Moreover,
$\Mod$ is locally factorial by (\ref{locally-factorial}) in this case.

\subsection{} Consider $G=SO_{r}$ with its standard (orthogonal) representation
and suppose that
$r\geq 7$. The moduli space $\Pic(\ModSO)$ is the good quotient of a
parameter scheme $\Quad$ by
$GL(H)$ with $H=k^{rN}$ (\cf \cite{So1}). The scheme $\Quad$ parameterizes
equivalent (with the obvious equivalence relation) triples
$([F,\sigma,\alpha])$, where $(F,\sigma)$ is a semistable
$SO_{r}$-bundle and
$\alpha:H^{0}(X,F(N))\isom H$.

Choose a theta-characteristic $\kappa$ on $X$. Then on $\Quad$ there is
the $GL(H)$-linearized pfaffian of cohomology line bundle ${\scr{P}}_{\kappa}$
deduced from the universal family over $\Quad\times X$.

\begin{prop}\label{Descent-of-P-to-Mreg} The line bundle ${\scr{P}}_{\kappa}$
descends to
$\ModSOreg$.
\end{prop}

\begin{proof} We use Kempf's lemma. If $r$ is even, the stabilizer at a point
$q=[F,\sigma,\alpha]\in\Quad^{reg}$ is $\pm 1$; if $r$ is odd, the stabilizer
is reduced to $1$. So in the latter case, there is nothing to prove. In the
former case, by definition of the pfaffian of cohomology, using that its
formation commutes with base change, the action $\pm 1$ is given by
$g\mapsto g^{h^{1}(F\otimes\kappa)}$, so the action is trivial, as
$h^{1}(F\otimes\kappa)$ is even.
\end{proof}

\begin{prop}\label{not-locally-factorial} If $r\geq 7$, the line bundle
${\scr{P}}_{\kappa}$ does not descend to
$\ModSOzero$. In particular, $\ModSOzero$ is not locally factorial.
\end{prop}

\begin{proof} Let $(F_{1},\sigma_{1})$ be a regularly stable {\em odd}
$SO_4$-bundle, and $(F_{2},\sigma_{2})$ be a regularly stable {\em odd}
$SO_{r-4}$-bundle. If $r=8$, suppose that $(F_{1},\sigma_{1})$ and
$(F_{2},\sigma_{2})$ are not isomorphic. Then the orthogonal sum
$(F,\tau)=(F_{1}\oplus F_{2},\sigma_{1}\oplus\sigma_{2})$ is {\em even}. Let
$[F,\tau,\alpha]\in\Quad$ be a point corresponding to $(F,\tau)$. Again,
by definition of the pfaffian of cohomology, using that its formation commutes
with base change,  we see that the action of the stabilizer $\{\pm
1\}\times\{\pm 1\}$  is
$$(g_{1},g_{2})\mapsto g_{1}^{h^{1}(F_{1}\otimes\kappa)}
g_{2}^{h^{1}(F_{2}\otimes\kappa)}.$$ But then the element $(-1,1)$ acts
nontrivially.
\end{proof}

\bigskip\bigskip

\yladress
\medbreak
\csadress

\end{document}